\documentclass[oneside,11pt]{article}

\usepackage[utf8]{inputenc}
\usepackage[T1]{fontenc}
\usepackage{amsthm}
\usepackage{amsmath}
\usepackage{amsfonts}
\usepackage{amssymb}
\DeclareMathOperator{\sech}{sech}
\usepackage{bbm}
\usepackage{enumitem}
\usepackage{graphicx} 
\usepackage{float}
\usepackage{booktabs}
\usepackage{multirow}
\usepackage{rotating}
\usepackage{array}
\usepackage{xcolor}
\usepackage[ruled]{algorithm2e}
\usepackage[caption=false]{subfig}
\newcolumntype{R}[1]{>{\raggedleft\let\newline\\\arraybackslash\hspace{0pt}}m{#1}}
\newcolumntype{L}[1]{>{\raggedright\let\newline\\\arraybackslash\hspace{0pt}}m{#1}}
\usepackage{breakcites}
\usepackage{nowidow}
\usepackage[top=1.5cm,bottom=2cm,right=2.5cm,left=2.5cm]{geometry}
\usepackage{titling}
\usepackage{natbib}
\usepackage{ifplatform}
\usepackage{textcomp}
\ifwindows
\fi
\allowdisplaybreaks

\newcommand{\R}{\mathbb{R}} 
\newcommand{\cF}{\mathcal{F}} 
\newcommand{\cD}{\mathcal{D}} 
\newcommand{\cC}{\mathcal{C}} 
\newcommand{\cB}{\mathcal{B}} 
\newcommand{\cA}{\mathcal{A}} 
\newcommand{\cR}{\mathcal{R}} 
\newcommand{\rmd}{\mathrm{d}} 
\newcommand{\InfGen}{\mathbb{L}} 
\newcommand{\Ind}{\mathbbm{1}} 
\newcommand{\lp}{\left(} 
\newcommand{\rp}{\right)} 
\newcommand{\lb}{\left[} 
\newcommand{\rb}{\right]} 
\newcommand{\lc}{\left\{} 
\newcommand{\rc}{\right\}} 
\newcommand{\blp}{\big(}
\newcommand{\brp}{\big)}
\newcommand{\blb}{\big[}
\newcommand{\brb}{\big]}
\newcommand{\blc}{\big\{}
\newcommand{\brc}{\big\}}
\newcommand{\lrav}[1]{\left|#1\right|} 
\newcommand{\wt}[1]{\widetilde{#1}}

\newcommand{\lrp}[1]{\lp#1\rp} 
\newcommand{\lrc}[1]{\lc#1\rc} 
\newcommand{\blrp}[1]{\blp#1\brp} 
\newcommand{\blrb}[1]{\blb#1\brb} 
\newcommand{\blrc}[1]{\blc#1\brc} 
\newcommand{\Prob}[1]{\mathbb{P}\lp #1\rp} 
\newcommand{\Pro}[2]{\mathbb{P}_{#2}\lp #1\rp} 
\renewcommand{\Pr}{\mathbb{P}} 
\newcommand{\Esp}[1]{\mathbb{E}\lb #1\rb} 
\newcommand{\Es}[2]{\mathbb{E}_{#2}\lb #1\rb} 
\newcommand{\E}{\mathbb{E}} 
\newcommand{\interior}[1]{%
  {\kern0pt#1}^{\mathrm{o}}%
}

\newtheorem{theorem}{Theorem}

\newtheorem{remark}{Remark}
\newtheorem{proposition}{Proposition}
\newtheorem{lemma}{Lemma}

\usepackage[pdftex,final,bookmarksnumbered,bookmarksopen=false,breaklinks,colorlinks]{hyperref}
\hypersetup{
	final,
	bookmarksnumbered=true,
	bookmarksopen=true,
	bookmarksopenlevel=0,
	unicode=false,
	pdftoolbar=true,
	pdfmenubar=true,
	pdffitwindow=false,
	pdftitle={Optimal exercise of American options under time-dependent Ornstein-Uhlenbeck processes},
	pdfdisplaydoctitle=true,
	pdftoolbar=true,
	pdfmenubar=true,
	pdflang={English},
	pdfauthor={Abel Azze, Bernardo D'Auria, Eduardo Garcia-Portugues},
	pdfsubject={arXiv paper},
	pdfcreator={Abel Azze},
	pdfproducer={Abel Azze},
	pdfkeywords={American option}{Free-boundary problem}{Optimal stopping}{Ornstein-Uhlenbeck}{Time-inhomogeneity},
	pdfnewwindow=true,
	breaklinks=true,
	hidelinks,       
	linkcolor=black,
	citecolor=black,
	filecolor=magenta,
	urlcolor=cyan
}

\newif\ifmain
\maintrue
\newif\ifsupplement
\supplementtrue

\newif\iffigstabs
\figstabstrue

\begin{document}

\ifmain

\title{Optimal exercise of American options under time-dependent Ornstein--Uhlenbeck processes}
\setlength{\droptitle}{-1cm}
\predate{}%
\postdate{}%
\date{}

\author{Abel Azze$^{1,4}$, Bernardo D'Auria$^{2}$, and Eduardo Garc\'ia-Portugu\'es$^{3}$}
\footnotetext[1]{Department of Quantitative Methods, CUNEF Universidad (Spain).}
\footnotetext[2]{Department of Mathematics ``Tullio Levi Civita'', University of Padova (Italy).}
\footnotetext[3]{Department of Statistics, Universidad Carlos III de Madrid (Spain).}
\footnotetext[4]{Corresponding author. e-mail: \href{mailto:abel.guada@cunef.edu}{abel.guada@cunef.edu}.}
\maketitle

\begin{abstract}
    We study the barrier that gives the optimal time to exercise an American option written on a time-dependent Ornstein--Uhlenbeck process, a diffusion often adopted by practitioners to model commodity prices and interest rates. By framing the optimal exercise of the American option as a problem of optimal stopping and relying on probabilistic arguments, we provide a non-linear Volterra-type integral equation characterizing the exercise boundary, develop a novel comparison argument to derive upper and lower bounds for such a boundary, and prove its Lipschitz continuity in any closed interval that excludes the expiration date and, thus, its differentiability almost everywhere. We implement a Picard iteration algorithm to solve the Volterra integral equation and show illustrative examples that shed light on the boundary's dependence on the process's drift and volatility.
\end{abstract}
\begin{flushleft}
\end{flushleft}

\section{Introduction}\label{sec:intro}

The extra flexibility of being able to exercise at any time before maturity makes American options a popular financial instrument among traders, but it makes these types of vanilla options significantly harder to price when compared to their European counterparts. Indeed, while the arbitrage-free price of the European option written on a geometric Brownian motion was settled down in the notorious works of \cite{Merton-1973-theory} and \cite{Black-1973-pricing}, the valuation of the American option had to go through a much longer process that took 40 years to complete, starting from the early works of \cite{Samuelson-1965-rational} and \cite{McKean-1965-free-boundary}, reaching important milestones in the works of \cite{Kim-1990-analytic}, \cite{Jacka-1991-optimal}, and \cite{Carr-1992-alternative}, and concluding in \cite{Peskir-2005-American}. This pricing journey is registered in the surveys of \cite{Myneni-1992-pricing} and \cite{Barone-Adesi-2005-saga}.

The classic geometric Brownian motion is not a suitable model in many financial contexts. Examples are the so-called mean-reverting strategies that arise when trading assets that tend to revert to an average price. \textit{Pair-trading} strategies are the most popular of these scenarios. They take place when the spread between correlated assets is pulled towards a baseline level. In such a case, it is common to model the spread by means of an Ornstein--Uhlenbeck (OU) process to capture the mean-reverting behavior. See \cite{Leung-2015-optimal-book} for theory and empirical evidence, \cite{Gatev-2006-pairs} for a practitioner perspective, and the handbook of \cite{Ehrman-handbook-2006} for a quick reference guide. Mean-reverting strategies also find room when trading options. Indeed, imagine that a trader holds two American options, one put and one call, written on two different stocks whose price at time $t$ is $X_t^{(1)}$ and $X_t^{(2)}$, with the same maturity date $T$ and strike prices given by $A_1$ and $A_2$. At any time $t$ before maturity, the trader can exercise both options, thus making the profit $\smash{\blrp{\alpha_1X_t^{(1)} - A_1} + \blrp{A_2 - \alpha_2X_t^{(2)}} = \lrp{A_2 - A_1} - \blrp{\alpha_2X_t^{(2)} - \alpha_1X_t^{(1)}}}$, where $\alpha_1$ and $\alpha_2$ represent the shares of each asset held by the trader in each American option position. The trader can then carve out a new American put option written on the spread $\smash{X_t = X_t^{(2)} - X_t^{(1)}}$ and with strike price given by $A = A_2 - A_1$, and the spread could be modeled as an OU process by conveniently choosing $\alpha_1$ and $\alpha_2$ (see Chapter 2 in \cite{Leung-2015-optimal-book}). 

An OU process has also been shown convenient for modeling certain commodity prices (see, e.g., \cite{Chaiyapo-2017-application} and \cite{Ogbogbo-2018-modeling}). However, it is far more common to use it to model log-prices (see, e.g., \cite{Schwartz-1997-stochastic}, \cite{Boyarchenko-2007-perpetual}, \cite{Zhang-2012-efficient}, and \cite{Mejia-2018-calibration}), which removes the drawback of allowing negative prices. Further flexibility can be added to the OU model if its parameters are time-dependent. Some authors, for instance, have relied on a time-dependent OU process to capture the effect of seasonality on commodity prices (see, e.g., \cite{Zapranis-2008-modelling}, \cite{Back-2013-seasonality}, and \cite{Tong-2021-censored}) and to cope with predictable assets \citep{Lo-1995-implementing, Carmona-2012-does}. A time-dependent OU process is also a go-to for practitioners when modeling interest rates, as it accommodates the popular Hull--White model \citep{Hull-1990-pricing}, an iteration of the Vasicek model (i.e., an OU process with constant parameters).

The analytical valuation of American options is fundamentally an Optimal Stopping Problem (OSP), as \cite{McKean-1965-free-boundary} established, and \cite{Bensoussan-1984-theory} and \cite{Karatzas-1988-pricing} later proved for the geometric Brownian motion case. The solution of the OSP is both the option's price, often called the value function, and the Optimal Stopping Boundary (OSB), that is, the optimal exercise barrier that sets the rule for deciding, at any time before maturity, whether it is optimal or not to exercise. Due to the connection between OSPs and free-boundary problems that started to be noticed in the 1960s by works like \cite{McKean-1965-free-boundary}, \cite{Grigelionis-1966-StefanProblem}, and \cite{Grigelionis-1967-optimal}, and later detailed, among others, in the book by \cite{Peskir-2006-optimal}, the valuation of an American option can also be seen through the lens of partial differential equations and, consequently, embedded into the well-consolidated theory of variational inequalities. The works by \cite{Jaillet-1990-variational} and \cite{Blanchet-2006-continuity} go in that direction. The downside of variational inequality methods, compared to probabilistic ones, is that they tend to lead to less explicit solutions of both the value function and the~OSB.

To bypass the complexity of pricing American options analytically, a plethora of numerical methods based on Monte Carlo simulation, finite differences, and binomial/trinomial trees have been developed. \cite{Zhao-2018-American} performs a comparison between eight of the most popular of these methods. Nevertheless, if available, analytical approaches are typically preferred over (pure) numerical ones, as the former are generally more precise and computationally expedient. Working out analytical solutions is not trivial, even if the American option is written on top of relatively simple models like the plain geometric Brownian motion. It becomes even harder when the coefficients of the underlying process are time-dependent. In this case the OSB cannot be guaranteed to be monotonic or convex, and proving the so-called \textit{smooth-fit condition} (the extra boundary condition that is required to offset not knowing the boundary in advance) becomes mathematically challenging. Under a geometric Brownian motion scheme, \cite{Ekstrom-2004-properties} worked with a time-dependent volatility and proved smoothing properties of the price function. \cite{Rehman-2009-American} took Ekstr\"om's work further by also allowing a time-dependent interest rate and providing the fair price of the American option in terms of the early exercise premium characterization, that is, the price of the European option plus the early exercise premium associated with the possibility of exercising before maturity. They did not prove the smooth-fit condition, which led to a somewhat difficult proof of the vanishing of the local-time term that the Itô--Tanaka formula yields for the pricing formula. Escaping the framework of a geometric Brownian motion, \cite[Section 4.6]{Detemple-2005-American-style} provided the early exercise premium characterization of the American option's price written on a diffusion with time- and space-dependent coefficients, these having somehow restrictive conditions that lead to a monotonic OSB. Using variational inequalities and partial differential equation techniques, \cite{Jaillet-1990-variational} and \cite{Blanchet-2006-continuity} also worked with time-inhomogeneous diffusions. 

In pricing American options, OU processes have been mostly used as a model for interest rates and discounting factors, rather than to model the underlying asset price. To consolidate their use in the latter case, we offer in this paper the solution to the problem of optimally exercising an American put option written on top of a time-inhomogeneous OU process with continuously differentiable coefficients. The associated OSB is given in terms of a non-linear Volterra-type integral equation that holds an explicit and relatively simple dependency on the OU parameters, which allows deriving a clean intuition on how those affect the OSB shape. We comment on two algorithms and implement one of them to work out the solution of the integral equation. The solution method we used, which is primarily probabilistic, yields important properties about the OSB, such as its Lipschitz continuity and, thus, differentiability almost everywhere. This order of smoothness of the OSB, which goes beyond its typical differentiability almost everywhere, is usually difficult to prove, especially for non-convex OSBs and non-differentiable gain functions. We also derive explicit upper and lower bounds for the OSB by relying on a comparison method that, to the best of our knowledge, has not been used within the framework of optimal stopping and could potentially be used in generic processes with drifts that are monotone in the space component. The value function is characterized by means of the early-exercise premium representation, which has a clear economic interpretation. 
The shape of the OSB is explored numerically for different sets of values of the drift and volatility of the OU process. We also cover the American call option scenario by proving that its OSB is the reflection of the OSB of its put counterpart with respect to the strike price.

Recently, by using a different approach, \cite{Carr-2021-semi-analytical} also tackled the problem of optimally exercising an American option with an OU process with differentiable time-dependent coefficients underneath. They performed time-space transformations to cast the partial differential operator of the associated free-boundary problem into that of the heat equation, which they solved by means of the heat potential method; the recent paper of \cite{Lipton-2020-physics} explains this method within the context of the American put option for a geometric Brownian motion. As a result, the OSB is proved to satisfy an integral equation of Fredholm type, although the uniqueness of the solution was not addressed, nor numerical methods for finding approximate solutions were discussed. In addition, the dependence of the integral equation with respect to the process's coefficients does not seem to be evident, as neither is the interpretation of the value formula in financial terms. Remarkably, \cite{Carr-2021-semi-analytical} also allow for a time-dependent discount rate and repulsion drifts, which we do not account for in our setting.

The rest of the paper is structured as follows. Section \ref{sec:formulation} sets the problem and introduces the basic notation. Section \ref{sec:properties} gathers all the properties required to come up, in Section \ref{sec:pricing_formula}, with the value formula and the free-boundary equation. Numerical experiments are displayed in Section \ref{sec:numerical_results}, while concluding remarks are relegated to Section \ref{sec:conclusions}.

\section{Formulation of the problem}\label{sec:formulation}
	
Define $\mu:[0,T]\times\mathbb{R}\rightarrow\mathbb{R}$ such that $\mu(s, x) := \theta(s)(\alpha(s) - x)$, where $\theta, \alpha:[0, T]\rightarrow \R$ are continuously differentiable functions, and $\theta(s) > 0$ for all $s\in[0, T]$. Let $X = (X_s)_{s \in [0, T]}$ be the unique strong solution (see, e.g., \citep[Theorem 2.1, Chapter IX]{Revuz-2010-stochastic}) of the stochastic differential equation
\begin{align}\label{eq:SDE}
    \rmd X_s = \mu(s, X_s)\,\rmd s + \rmd W_s,\ \ 0 \leq s \leq T,
\end{align}
in the filtered space $(\Omega, \cF, \Pr, (\cF_s)_{s \in [0, T]})$, where $(\cF_s)_{s \in [0, T]}$ is the natural filtration of the underlying standard Brownian motion $(W_s)_{s \in [0, T]}$. We refer to $\theta$ and $\alpha$ as the slope and pulling level (functions), respectively.

From \eqref{eq:SDE}, we get that the infinitesimal generator $\InfGen$ of the process $((s, X_s))_{s \in [0, T]}$ is given by
\begin{align}\label{eq:OUB_InfGen}
    (\InfGen f)(s, x) = \partial_t f(s, x) + \theta(s)(\alpha(s) - x)\partial_x f(s, x) + \frac{1}{2} \partial_{xx}f(s, x).
\end{align}
Above and henceforth, we use $\partial_t$ and $\partial_x$ to denote, respectively, the usual partial derivatives with respect to time and space of a function with arguments $(t,x)$, and $\partial_{xx}$ is used as a shorthand for $\partial_x\partial_x$.

Consider the finite-horizon, discounted Optimal Stopping Problem (OSP)
\begin{align}\label{eq:OSP}
    V(t, x) = \sup_{0\leq \tau \leq T - t}\Es{e^{-\lambda \tau}G(X_{t + \tau})}{t, x},
\end{align}
where $V$ is the value function, 
\begin{align}\label{eq:gain_function}
    G(x) := (A - x)^+
\end{align}
is the gain function for some $A\in\R$ that represents the strike price of the option, and $\lambda \geq 0$ is the discounting rate. The supremum above is taken over all random times $\tau$ such that $t + \tau$ is a stopping time in $(\cF_s)_{s \in [0, T]}$ and $\E_{t, x}$ represents the expectation under the probability measure $\Pr_{t, x}$ defined as $\Pro{\cdot}{t, x} := \Prob{\cdot\ |\ X_t = x}$. For simplicity and in an abuse of notation, in the sequel, we will refer to $\tau$ as a stopping time while keeping in mind that $t + \tau$ is the actual stopping time in the filtration $(\cF_s)_{s \in [0, T]}$. The particular form of $G$ implies that
\begin{align}\label{eq:G_inequality}
    -(y - x)^- \leq G(x) - G(y) \leq (y - x)^+ 
\end{align}
for all $x, y \in \R$, inequalities which will be recurrently used throughout the paper.

It is useful to keep track of the condition $X_t = x$ in a way that does not change the probability measure whenever $t$ or $x$ changes. To do so, we denote the process \mbox{$X^{t, x} = (X_s^{t, x})_{s \in [0, T-t]}$} in the filtered space $(\Omega, \cF, \Pr, (\cF_s)_{s \in [0, T-t]})$ as the unique solution~of
\begin{equation*}
    \left\{
    \begin{aligned}
    \rmd X_s^{t, x} &= \mu(t + s, X_s^{t, x})\,\rmd s + \rmd W_s,\ \ 0 \leq s \leq T - t, \\
    X_0^{t, x} &= x,
    \end{aligned}
    \right.
\end{equation*}
where the underlying Brownian motion $\lrp{W_s}_{s\in\R_+}$ remains the same for all $t\in[0, T]$. Theorem~38 from Chapter~V in \cite{Protter-2003-stochastic} guarantees that $X^{t, x}$ is well defined and continuous with respect to $x$.

Theorem 39 from Chapter V in \cite{Protter-2003-stochastic} ensures that $X^{t, x}$ is continuously differentiable with respect to $x$ and $t$, and combined with \eqref{eq:SDE}, it states that the processes \mbox{$\partial_t X^{t, x} = (\partial_t X_s^{t, x})_{s \in [0, T-t]}$} and $\partial_x X^{t, x} = (\partial_x X_s^{t, x})_{s \in [0, T-t]}$, defined as the following $\Pr$-a.s. limits
\begin{align*}
    \partial_t X_s^{t, x} := \lim_{\varepsilon\rightarrow 0} \left(X_s^{t + \varepsilon, x} - X_s^{t, x}\right)\varepsilon^{-1} \quad\text{and}\quad
    \partial_x X_s^{t, x} := \lim_{\varepsilon\rightarrow 0} \left(X_s^{t, x + \varepsilon} - X_s^{t, x}\right)\varepsilon^{-1},
\end{align*}
take the form
\begin{align}
    \partial_t X_s^{t, x} &= \int_{0}^{s}\left(\partial_t\mu(t + u, X_u^{t, x}) + \partial_x\mu(t + u, X_u^{t, x})\partial_t X_u^{t, x}\right)\,\rmd u \nonumber \\
    &= \int_{0}^{s}\lrp{\theta'(t + u)(\alpha(t + u) - X_u^{t, x}) + \theta(t + u)(\alpha'(t + u) - \partial_t X_u^{t, x})}\,\rmd u \label{eq:SDE_t} 
\end{align}
and
\begin{align}
    \partial_x X_s^{t, x} &= 1 + \int_{0}^{s}\partial_x\mu(t + u, X_u^{t, x})\partial_x X_u^{t, x}\,\rmd u \nonumber  \\
    &= 1 - \int_{0}^{s}\theta(t + u)\partial_x X_u^{t, x}\,\rmd u. \label{eq:SDE_x}
\end{align}
Mind the difference between the differential operators $\partial_t$ and $\partial_x$ and the processes $\partial_t X_s^{t, x}$ and $\partial_x X_s^{t, x}$. From \eqref{eq:SDE_x}, it follows that the process $\partial_x X_s^{t, x}$ is deterministic, i.e.
\begin{align}\label{eq:SDE_x^}
    \partial_x X_s^{t, x} &= \exp\left\{-\int_{0}^{s}\theta(t + u)\,\rmd u\right\}, 
\end{align}
and that $\partial_x X_s^{t, x}\in(0,1]$ for all $s\in[0,T-t]$.

The results in Section 2.2 from \cite{Peskir-2006-optimal}, along with the fact that $X$ is a Feller process (see, e.g., \cite[Theorem 2.5, Chapter IX]{Revuz-2010-stochastic}), give us the existence and characterization of a stopping time that is optimal in \eqref{eq:OSP}. Specifically, if we denote by $\cD := \{V = G\}$ and $\cC := \cD^c =\{V > G\}$ to the so-called stopping set and continuation set, respectively, then, if 
\begin{align}\label{eq:Dominated_convergence_condition}
    \E_{t, x}\bigg[\sup_{0\leq s\leq T - t} e^{-\lambda s}G(X_{t + s})\bigg] < \infty
\end{align} 
for all $(t, x) \in [0, T)\times\R$, we can guarantee that, under $\Pr_{t, x}$, the first hitting time of $(X_{t + s})_{s \in [0, T-t]}$ into $\cD$, denoted by $\tau^* = \tau^*(t, x)$, satisfies
\begin{align}\label{eq:OST}
    V(t, x) = \Es{e^{-\lambda\tau^*}G(X_{t + \tau^*})}{t, x}.
\end{align}
If there is another Optimal Stopping Time (OST) $\tau$, then $\tau^* \leq \tau\ \ \Pr_{t, x}$-a.s. (see Equation 2.2.1 from \cite{Peskir-2006-optimal}). Solving the OSP \eqref{eq:OSP} means providing tractable expressions for both the value function $V$ and the OST $\tau^*$.

The boundary of $\cD$ (or $\cC$), denoted by $\partial \cD$ (or $\partial \cC$), is called the Optimal Stopping Boundary (OSB). It turns out that both the value function $V$ and the OSB $\partial \cD$ are the solution of a Stefan problem with the infinitesimal generator $\InfGen$ being the differential operator acting on $V$ in $\cC$. For this reason, the OSB is also referred to as the free boundary; see \cite{Peskir-2006-optimal} for more information on the relation between OSPs and free-boundary problems. If the OSB can be depicted by the graph of a function $b:[0, T] \rightarrow \R$, i.e., if $\partial \cD = \{(t, b(t)): t \in [0, T]\}$, then $b$ is referred to as the OSB too.

Finally, it is convenient to recall the martingale and supermartingale properties of~$V$:
\begin{align}
    \Es{V(t + \tau^*\wedge s, X_{t + \tau^*\wedge s})}{t, x} &=  V(t, x), \label{eq:V_martingale} \\
    \Es{V(t + s, X_{t + s})}{t, x} &\leq V(t, x), \label{eq:V_supermartingale}
\end{align}	
for all $0 \leq s < T - t$, as they will be used to prove results in Section \ref{sec:properties}.

\section{Regularities of the boundary and the value function}\label{sec:properties}

In this section we derive regularity conditions about the OSB and the value function. These conditions allow obtaining a solution for the OSP \eqref{eq:OSP}, later addressed in Section~\ref{sec:pricing_formula}, by using an extension of the Itô formula to derive a characterization of the OSB via a Volterra integral equation. 

We first point out that \eqref{eq:Dominated_convergence_condition} holds true, as it is a consequence of Lemma \ref{lm:bounds} below. This allows us to prove the existence of the OST as well as its characterization in terms of the stopping set claimed in Section \ref{sec:formulation}. We then shed light on the geometry of the stopping and continuation regions in the next proposition, which entails that the stopping set lies below the continuation set and that the boundary between them can be seen as the graph of a bounded function. We highlight the comparison argument used to derive the lower bound of the OSB. Besides its simplicity, it relies only on the monotonicity of the drift with respect to the space component, which makes it usable beyond the specifics of the OU process.

\begin{proposition}[Boundary existence and shape of the stopping set]\label{pr:boundary_existence}\ \\
	There exists a function $b:[0, T]\rightarrow\mathbb{R}$ such that $-\infty < b(t) < A$ for all $t\in[0, T)$ and $\cD = \{(t, x)\in[0, T)\times\mathbb{R} : x \leq b(t)\}\cup \lrp{\lrc{T}\times\R}$, where $A$ is as in \eqref{eq:gain_function}. Moreover, \mbox{$b(t) \leq \gamma(t)$} for all $t\in[0, T)$, where 
    \begin{align}\label{eq:gamma}
        \gamma(t) := (\theta(t) + \lambda)^{-1}(\theta(t)\alpha(t) + \lambda A)
    \end{align}
    and $b(T) := \lim_{t\rightarrow T} b(t) = \min\{A, \gamma(T)\}$. 
\end{proposition}

\begin{proof}
	Define $b:[0, T]\rightarrow\R$ such that $b(t) := \sup\{x \in\mathbb{R}: (t, x) \in \cD\}$ for all $t\in[0, T)$ and $b(T) = \lim_{t\rightarrow T}b(t)$. We first prove that $b(t) < A$ for all $t\in[0,T)$. Let $(t, x)\in \cA := [0, T)\times[A, \infty)$ and define the stopping time $\tau_\delta := \inf\{s\in[0, T - t] :\ X_{t + s} \leq A - \delta\}$ for some $\delta > 0$, assuming that $\inf\{\emptyset\} = T - t$. Since $\tau_{\delta} \leq T - t$ and $\Prob{\tau_\delta < T - t} > 0$,
	\begin{align*}
	    V(t, x) \geq \mathbb{E}_{t, x}\left[e^{-\lambda\tau_{\delta}}G(X_{t + \tau_{\delta}})\right] \geq e^{-\lambda(T - t)}\Prob{\tau_\delta < T - t}\delta > 0 = G(x),
	\end{align*}
	that is, $(t, x)\in \cC$. Therefore, $\cA \subset \cC$.
	
	We now prove that $\cD = \{(t, x)\in[0, T]\times\mathbb{R} : x \leq b(t)\}\cup \lrp{\lrc{T}\times\R}$. The fact that $\lrc{T}\times\R \subset \cD$ is a straightforward consequence of the definition $V(T, x) = G(x)$ for all $x\in\R$.
	Moreover, if $(t, x_1)\in \cD$, then, for any $x_2 < x_1$ we have that, for $\tau^* = \tau^*(t, x_2)$,
	\begin{align}
	    V(t, x_2) - V(t, x_1) &\leq \Esp{e^{-\lambda\tau^*}\lrp{G\lrp{X_{\tau^*}^{t, x_2}} - G\lrp{X_{\tau^*}^{t, x_1}}}} \nonumber \\
	    &\leq \Esp{X_{\tau^*}^{t, x_1} - X_{\tau^*}^{t, x_2}} 
	    = (x_1 - x_2) \Esp{\partial_x X_{\tau^*}^{t, x}} \label{eq:V_x1-V_x2} \\
	    &\leq x_1 - x_2 = G(x_2) - G(x_1), \nonumber
	\end{align}
        where in the second inequality we used \eqref{eq:G_inequality} and the fact that $x\mapsto X_s^{t, x}$ is an increasing function for all $s\in [0, T - t)$ ($\partial_x X_s^{t, x} \geq 0$ for all $s$ according to \eqref{eq:SDE_x^}). We also used the mean value theorem alongside \eqref{eq:SDE_x^} to come up with the equality.
	Since $V(t, x_1) = G(x_1)$, it follows from \eqref{eq:V_x1-V_x2} that $(t, x_2)$ belongs to $\cD$. Finally, being $\cD$ closed, $(t, b(t)) \in \cD$ for all $t\in[0, T)$, which guarantees that $\cD$ has the claimed shape.

	We now use the Itô--Tanaka formula to get that, for all $(t, x)\in[0, T)\times\R$,
    \begin{align}\label{eq:G_ItoTanaka}
        e^{-\lambda s}G(X_s^{t, x}) =&\; G(x) + \int_{0}^{s}F_t\lrp{u, X_u^{t, x}}\Ind(X_u^{t, x} < A)\,\rmd u  \\
        & - \int_{0}^{s}e^{-\lambda u}\Ind(X_u^{t, x} < A)\,\rmd B_{u} + \frac{1}{2}\int_0^{s}e^{-\lambda u}\Ind(X_u^{t, x} = A)\,\rmd l_u^A(X^{t, x}), \nonumber
    \end{align}
    where $F_t(u, x) = -e^{-\lambda u}\lrp{\lambda G(x) + \mu(t+u, x)}$, and $l_s^A(X)$ is the local time of the process $X$ at $A$ up to time $s$, that is,
	\begin{align*}
	l_s^A(X) = \lim_{h\downarrow 0} \frac{1}{2h}\int_{0}^{s}\Ind(A - h \leq X_u \leq A + h)\,\rmd \langle X, X\rangle_u.
	\end{align*}
    We have already proved that $b(t) \leq \gamma(t)$ whenever $\gamma(t) \geq A$, so we will assume that $\gamma(t) < A$. Fix~$(t, x)$ such that $\gamma(t)< x < A$ and notice that, in such a case,
    \begin{align*}
        F_t(u, x) > 0& \;\text{ if }\; x > \gamma(t+u), \\
        F_t(u, x) < 0& \;\text{ if }\; x < \gamma(t+u).
    \end{align*}
    Take a ball $\cB \subset [0, T]\times(-\infty, A)$ centered at $(t, x)$ and sufficiently small such that $F_t(u, y) > 0$ for all $(u, y)\in\cB$. In \eqref{eq:G_ItoTanaka}, replacing $s$ by the first exit time from $\cB$ (denoted by $\tau_{\cB^c}$), taking $\Pr_{t, x}$-expectation (cancels the martingale term), and noticing that $X_{t + u}$ never touches $A$ for $u\leq\tau_{\cB^c}$ under $\Pr_{t, x}$ (i.e., the local-time term is null), we get
    \begin{align*}
        V(t, x) \geq \Es{e^{-\lambda \tau_{\cB^c}}G(X_{t+\tau_{\cB^c}})}{t, x} = G(x) + \Es{\int_{t}^{t + \tau_{\cB^c}}F_t\lrp{u, X_u^{t, x}}\,\rmd u}{t, x} > G(x),
    \end{align*}
    meaning that $(t, x)\in\cC$ and, therefore, that $b(t) \leq \gamma(t)$ for all $t\in[0, T)$. 

    We now prove that $(t, x) \in \cD$ for all $x < \min\lrc{A, \gamma(T)}$ and all $t$ sufficiently close to $T$. This, along with the fact that $b(t) \leq \min\lrc{A, \gamma(t)}$ that we have already proved, shows that $\lim_{t\rightarrow T}b(t)$ exists and takes the claimed value. Notice that
    \begin{align}
        \Es{\int_{0}^{s}F_t\lrp{u, X_u^{t, x}}\Ind(X_u^{t, x} < A)\,\rmd u}{t, x} &\leq s\Es{\sup_{0\leq u\leq s}F_t\lrp{u, X_u^{t, x}}\Ind(X_u^{t, x} < A)}{t, x} \nonumber \\
        &\leq s \lrp{K_1 p(t, x, s, \varepsilon) + K_2 (1-p(t, x, s, \varepsilon))}\label{eq:G_ItoTanaka_1}
    \end{align}
    and
    \begin{align}\label{eq:G_ItoTanaka_2}
        \Es{\int_0^{s}e^{-\lambda u}\Ind(X_u^{t, x} = A)\,\rmd l_u^A(X^{t, x})}{t, x} \leq s(1-p(t, x, s, \varepsilon)),
    \end{align}
    with 
    \begin{align*}
        p(t, x, s, \varepsilon) := \Prob{\sup_{0\leq u\leq s} X_u^{t, x} < x+\varepsilon}
    \end{align*}
    and
    \begin{align*}
        K_1 := \sup_{0\leq u\leq s}F_t\lrp{u, x+\varepsilon} < 0, \quad K_2 := \lrp{\sup_{0\leq u\leq s}F_t\lrp{u, A}}^+,
    \end{align*}
    for $\varepsilon > 0$ and $t\in[0, T)$ such that $x + \varepsilon \leq \min \lrc{A, \min_{u\in[0, T-t]} \gamma(t+u)}$. This selection of $\varepsilon$ and $t$ is possible, as the continuity of $\gamma$, inherited from the same level of smoothness of $\theta$ and $\alpha$, guarantees that $\gamma(T)$ gets arbitrarily close to $\gamma(T)$ as $t$ converges to $T$. Likewise, $t$ close enough to $T$ yields $p(t, x, s, \varepsilon) \geq (K_2+1)/(K_2+1-K_1)$. Indeed,
    \begin{align*}
        \lim_{t\rightarrow T}\Prob{\sup_{0\leq u\leq s} X_u^{t, x} < x+\varepsilon} &\geq  \lim_{t\rightarrow T}\Prob{\sup_{0\leq u\leq T-t} X_u^{t, x} - x - \varepsilon < 0} \\
        & = \lim_{t\rightarrow T}\Prob{\sup_{0\leq u\leq T-t} \frac{X_u^{t, x} - x - \varepsilon}{2u\ln(\ln(1/u))} < 0} \\
        &= \Prob{\lim_{t\rightarrow T}\sup_{0\leq u\leq T-t} \frac{X_u^{t, x} - x - \varepsilon}{2u\ln(\ln(1/u))} < 0}\\
        &= \Prob{\limsup_{u\rightarrow0} \frac{X_u^{t, x} - x - \varepsilon}{2u\ln(\ln(1/u))} < 0} = 1,
    \end{align*}
    where the last inequality comes after \eqref{eq:OU-LIL_sup} below. Therefore, after using \eqref{eq:G_ItoTanaka} along with \eqref{eq:G_ItoTanaka_1} and \eqref{eq:G_ItoTanaka_2}, we get that
    \begin{align*}
        V(t, x) - G(x) &\leq \sup_{s\in[0, T-t]} s \lrp{K_1 p(t, x, s, \varepsilon) + (K_2+1) (1-p(t, x, s, \varepsilon))} \leq 0,
    \end{align*}
    meaning that $(t, x)\in \cD$.
    
    The lower boundedness of $b$ is now addressed. 
    Define $m := \inf_{t\in[0,T]}\alpha(t)$. The following inequalities hold $\Pr$-a.s.:
    \begin{align}\label{eq:X>B}
        X_s^{t, x} &\geq Z_s^{t, x} \geq \alpha(t+s) - \left|\alpha(t+s) - B_s^x\right| \geq m - \left|m - B_s^x\right|,
    \end{align}  
    where $Z_s^{t, x} := \alpha(t + s) - \big|\alpha(t + s) - X_s^{t, x}\big|$, $B_t^x = x + W_t$, and $(W_t)_{t\geq0}$ is the underlying standard Brownian motion in \eqref{eq:SDE}. The second inequality holds since the drift of $Z := (Z_s^{t, x})_{s\in[0, T]}$ is greater than the drift of the reflection of $B^x$ with respect to $\alpha$ for all $(t, x)\in[0, T)\times\R$, and therefore we can assume that the first process is greater than the last one pathwise $\Pr$-a.s. \citep[see][Corollary 3.1]{Peng-2006-necessary}. 
    Indeed, for $\varepsilon > 0$, define the function
    \begin{align*}
        g_\varepsilon(x,\alpha) := 
        \begin{cases}
        |x-\alpha|, & \text{ if } |x-\alpha| \geq \varepsilon, \\
        \frac{1}{2}\lrp{\varepsilon + \varepsilon^{-1}(x-\alpha)^2}, & \text{ if } |x-\alpha| < \varepsilon
        \end{cases}
    \end{align*}
    and the processes
    \begin{align*}
        Y_t^{(1), \varepsilon} := \alpha(t) - g_\varepsilon\lrp{X_t, \alpha(t)}, \quad Y_t^{(2),\varepsilon} := \alpha(t) - g_\varepsilon\lrp{W_t, \alpha(t)}.
    \end{align*}
    Denote by $\mu^{(i),\varepsilon}$ the drift of $Y^{(i), \varepsilon}$, $i = 1, 2$. By using the Itô formula, we obtain that
    \begin{align*}
        (\mu^{(1),\varepsilon} - \mu^{(2),\varepsilon})(t, x) = 
        \begin{cases}
        \theta(t)|\alpha(t) - x|, & \text{ if } |\alpha(t) - x| \geq \varepsilon, \\
        \varepsilon^{-1}\theta(t)(\alpha(t) - x)^2, & \text{ if } |\alpha(t) - x| < \varepsilon.
        \end{cases}
    \end{align*}
    Therefore, $\mu^{(1),\varepsilon} \geq \mu^{(2),\varepsilon}$, which allows us to use Corollary 3.1 from \cite{Peng-2006-necessary} to state that $Y_t^{(1), \varepsilon} \geq Y_t^{(2), \varepsilon}$ for all $t\in[0, T]$ $\Pr$-a.s. and obtain the claimed result after taking $\varepsilon \rightarrow 0$ and realizing that, in such a case, $g_\varepsilon(c, x) \downarrow |c - x|$.
    Indeed,
    \begin{align*}
        &\Prob{\sup_{s\in[0,T-t]} \lrc{Z_s^{t,x} - \lrp{\alpha(t+s) - \left|\alpha(t+s) - B_s^x\right|}} < 0} \\
        &\qquad= \Pro{\lim_{n\rightarrow \infty} \sup_{s\in[0,T-t]} \lrc{Y_{t+s}^{(1),1/n} - Y_{t+s}^{(2),1/n}} < 0}{t, x} \\
        &\qquad= \Pro{\bigcup_{N\geq 0} \bigcap_{n\geq N} \lrc{\sup_{s\in[0,T-t]} \lrc{Y_{t+s}^{(1),1/n} - Y_{t+s}^{(2),1/n}} < 0}}{t, x} = 0.
    \end{align*}
    In the first equality, the interchange of the limit and supremum operators is justified since the convergence of $g_\varepsilon$ is uniform, while the last equality comes since we are considering the probability of a numerable union and intersection of $\Pr_{t, x}$-null sets.
    
	The inequalities in \eqref{eq:X>B}, and the fact that $G$ is a decreasing function, guarantee that the value function \eqref{eq:OSP} is lower than the value function associated with the infinite-horizon, discounted OSP of a reflected (with respect to $m$) Brownian motion. Thus, the respective OSBs hold the reverse inequality, that is, the stopping region of the latter process is contained in the stopping region of the former. It is easy to show that the free boundary for the infinite-horizon, discounted OSP, with the gain function $G$ and a $m$-reflected Brownian motion, is finite. Actually, one can explicitly obtain the (constant) OSB by directly solving the associated free-boundary problem.
	Therefore $b$ is bounded from below.
\end{proof}

The value function satisfies the regularity properties
listed in the next proposition. The method used to obtain the Lipschitz continuity of $V$, based on properties 
\eqref{eq:V_martingale} and \eqref{eq:V_supermartingale}, 
is a powerful technique to solve OSPs framed within Markovian processes. For example, the work of \cite{DeAngelis-2019-Lipschitz} provides Lipschitz continuity of the value function in a similar fashion, but for differentiable gain functions.

\begin{proposition}[Regularity of $V$]\label{pr:V_regularity}\ \\
	The value function $V$ satisfies the following properties:
	\begin{enumerate}[label=(\roman{*}), ref=(\textit{\roman{*}})]
		\item \label{pr:V_Lcontinuous} $V$ is locally Lipschitz continuous.
		\item \label{pr:V_C^12onC} $V$ is $C^{1,2}$ on $\cC$ and on $\cD$, and $\mathbb{L}V = \lambda V$ on $\cC$. 
		\item \label{pr:x->V} $x\mapsto V(t, x)$ is decreasing for all $t\in [0, T)$. Moreover, for $x\neq b(t)$,
		\begin{align}\label{eq:V_x}
		0 \geq \partial_xV(t, x) \geq -\E\left[e^{-\lambda\tau^*}\partial_x X_{\tau^*}^{t, x}\right],
		\end{align} 
		where $\tau^* = \tau^*(t, x)$. Additionally,
		\begin{align}
			\partial_t V(t, x) &\leq \Esp{\big|\partial_tX_{\tau^*}^{t, x}\big|}, \label{eq:V_t<} \\
			\partial_t V(t, x) &\geq -\Esp{\big|\partial_tX_{\tau^*}^{t, x}\big|} - L\Pr\left(\tau^* = T - t\right), \label{eq:V_t>}
		\end{align}	
		for some positive constant $L$.
	\end{enumerate}
\end{proposition}

\begin{proof}
    \ref{pr:V_Lcontinuous}\ \ 
    Fix $(t, x)\in [0, T]\times \R$. Since $\tau^* = \tau^*(t, x)$ is optimal for the OSP with initial conditions~$(t, x)$ and suboptimal for any different initial condition, and due to inequality \eqref{eq:G_inequality}, it follows that, for $\delta > 0$ small enough,  
	\begin{align}
		V(t, x) - V(t - \delta, x)
		&\leq \E\left[G\left(X_{\tau^*}^{t, x}\right) - G\left(X_{\tau^*}^{t - \delta, x}\right)\right] \leq \E\left[\left(X_{\tau^*}^{t - \delta, x} - X_{\tau^*}^{t, x}\right)^+\right] \nonumber \\
		&= \delta\E\left[\left(-\partial_tX_{\tau^*}^{t_\delta, x}\right)^+\right] \leq \delta L_{x}, \label{eq:Vt-Vt-delta<}
	\end{align}
	where the last equality holds for a certain (random) time $t_\delta\in(t - \delta, t)$ due to the mean value theorem. In \eqref{eq:Vt-Vt-delta<} and thereafter we use the notation
	\begin{align*}
	    L_{x} := \sup_{t \in [0,T]}\E\bigg[\sup_{s \leq T - t}\left|\partial_tX_s^{t, x}\right|\bigg].    
	\end{align*}
	By \eqref{eq:partial_t_X<} below, the value of $L_x$ is bounded.
	Applying similar arguments as those used in \eqref{eq:Vt-Vt-delta<}, relying on the martingale and supermartingale properties of the value function, and noticing that $\tau^*\wedge(T - t - \delta)$ is admissible for $V(t - \delta, x)$, we have
	\begin{align}
	    V(t + \delta,&\; x) - V(t, x) \nonumber\\
		\geq&\; \Esp{V\left(t + \delta + \tau^*\wedge(T - t - \delta), X_{\tau^*\wedge(T - t - \delta)}^{t + \delta, x}\right)} - \Esp{V\left(t + \tau^*, X_{\tau^*}^{t, x}\right)} \nonumber \\
		=&\; \Esp{\Ind(\tau^* \leq T - t - \delta)\left(V\left(t + \delta + \tau^*, X_{\tau^*}^{t + \delta, x}\right) - V\left(t + \tau^*, X_{\tau^*}^{t, x}\right)\right)} \nonumber \\
		& +\Esp{\Ind(\tau^* > T - t - \delta)\left(V\left(T, X_{T - t - \delta}^{t + \delta, x}\right) - V\left(t + \tau^*, X_{\tau^*}^{t, x}\right)\right)} \nonumber \\
		\geq&\; \Esp{\Ind(\tau^*\leq T - t - \delta)\left(G\left(X_{\tau^*}^{t + \delta, x}\right) - G\left(X_{\tau^*}^{t, x}\right)\right)} \nonumber \\
		& +\Esp{\Ind(\tau^* > T - t - \delta)\lrp{G(X_{T - t - \delta}^{t + \delta, x}) - G(X_{T - t - \delta}^{t, x})}} \nonumber \\
		& +\Esp{\Ind(\tau^* > T - t - \delta)\lrp{G(X_{T - t - \delta}^{t, x}) - V(t + \tau^*, X_{\tau^*}^{t, x})}} \nonumber \\
		\geq& -\Esp{\Ind(\tau^*\leq T - t - \delta)\left(X_{\tau^*}^{t, x} - X_{\tau^*}^{t + \delta, x}\right)^-} \nonumber \\
		& - \Esp{\Ind(\tau^* > T - t - \delta)\left(X_{T - t - \delta}^{t, x} - X_{T - t - \delta}^{t + \delta, x}\right)^-}\nonumber \\
		& +\Esp{\Ind(\tau^* > T - t - \delta)\left(G(X_{T - t - \delta}^{t, x}) - \Esp{V\lrp{T-\delta + \rho^*, X_{\rho^*}^{T-\delta, X_{T - t - \delta}^{t, x}}}}\right)} \nonumber \\
		\geq& -\delta L_{x} + \Esp{\Ind(\tau^* > T - t - \delta)\left(G(X_{T - t - \delta}^{t, x}) - V(T - \delta, X_{T - t - \delta}^{t, x})\right)}, \label{eq:Vt+delta-Vt>}
	\end{align}
	where $\rho^* := \tau^*\blrp{T - \delta, X_{T - t - \delta}^{t, x}}$ and, in the second inequality, we added and subtracted the term $G(X_{T - t - \delta}^{t, x})$. For the third inequality, we used the tower property of the expectation along with the Markovian nature of $X$, to get
    \begin{align*}
        \Esp{\Ind(\tau^* > T - t - \delta)\lrp{V(t + \tau^*, X_{\tau^*}^{t, x}) - \Esp{V\lrp{T-\delta + \rho^*, X_{\rho^*}^{T-\delta, X_{T - t - \delta}^{t, x}}}}}} = 0.
    \end{align*}
    We now analyze the last term in \eqref{eq:Vt+delta-Vt>} in detail. By the definition of $V$ in \eqref{eq:OSP}, using the Itô--Tanaka formula, and acknowledging that $\big(T - \delta, X_{T - t - \delta}^{t, x}\big) \in \cC$ in the set $\{\tau^* > T - t - \delta\}$, we derive the following inequality for $\rho^*$ in $\{\tau^* > T - t - \delta\}$ and~$Y_u = X_u^{T - \delta, X_{T - t - \delta}^{t, x}}$:
	\begin{align}
		V(T - \delta,&\; X_{T - t - \delta}^{t, x}) - G(X_{T - t - \delta}^{t, x}) \nonumber \\
		=&\; \Esp{e^{-\lambda\rho^*}G\lrp{Y_{\rho^*}}} - G(X_{T - t - \delta}^{t, x})\nonumber \\
		=&\; \mathbb{E}\Bigg[\frac{1}{2}\int_{0}^{\rho^*}e^{-\lambda u}\Ind\lrp{Y_u = A}\,\rmd l_u^A\lrp{Y}\nonumber\\
		&- \int_{0}^{\rho^*}e^{-\lambda u}\Ind\lrp{Y_u < A}\lrp{\lambda\lrp{A - Y_u} + \mu\lrp{T - \delta + u, Y_u}}\,\rmd u\Bigg] \nonumber \\
		\leq&\; \Esp{\frac{1}{2}\int_{0}^{\rho^*}\Ind\lrp{Y_u = A}\,\rmd l_u^A\lrp{Y} - \int_{0}^{\rho^*}e^{-\lambda u}\Ind\lrp{Y_u < A}\mu\lrp{T - \delta + u, Y_u}\,\rmd u} \nonumber \\
		\leq&\; \delta L, \label{eq:V_Lcont}
	\end{align}
	where $L = \frac{1}{2} + \max\{|\mu(t, x)| : 0\leq t \leq T,\ b(t) \leq x \leq A\} <  \infty$.
	Plugging \eqref{eq:V_Lcont} into \eqref{eq:Vt+delta-Vt>}, we obtain
	\begin{align}
	V(t + \delta, x) - V(t, x)	\geq - \delta \lrp{L_{x} + L\Pr\left(\tau^* > T - t - \delta\right)}. \label{eq:Vt+delta-Vt>_1}
	\end{align} 
	
	Now consider $\tau_\delta = \tau^*(t - \delta, x)$ for $0 \leq \delta \leq t \leq T$, and notice that we get the following from \eqref{eq:Vt+delta-Vt>_1}:
    \begin{align*}
        V(t, x) - V(t - \delta, x) \geq - \delta \lrp{L_{x} + L}.
    \end{align*}
	
    For $0 \leq \delta \leq T - t $ and $\tau^\delta = \tau^*(t + \delta, x)$ one gets the following by proceeding as in~\eqref{eq:Vt-Vt-delta<}:
	\begin{align*}
		V(t + \delta, x) - V(t, x) &\leq  \delta L_x.
	\end{align*}
			
	So far, since $L_{x}$ is finite (see Lemma \ref{lm:bounds}), we have proved that, for any $x\in\R$, $t\mapsto V(t, x)$ is Lipschitz continuous on $[0, T]$. We will now prove that $V$ is also Lipschitz continuous with respect to $x\in \R$ for all $t\in [0,T]$, which will complete the proof of~\ref{pr:V_Lcontinuous}.
			
	Since $G$ is decreasing and $x\mapsto X_s^{t, x}$ is increasing for all $s\in[0, T - t)$ and $t\in[0, T)$ ($\partial_x X_s^{t, x} \geq 0$), then $x\mapsto V(t, x)$ is decreasing for all $t\in[0, T)$. Fix $(t, x)\in[0, T)\times\R$ and $\delta > 0$. Consider $\tau^* = \tau^*(t, x)$, recall \eqref{eq:SDE_x^}, and argue as in \eqref{eq:V_x1-V_x2} to get
	\begin{align}
		0 \geq V(t, x + \delta) - V(t, x) &\geq -\delta\E\left[e^{-\lambda\tau^*}\exp\lrc{-\int_0^{\tau^*}\theta(t+u)\,\rmd u}\right] \label{eq:Vx+delta-Vx}
		\geq -\delta. 
	\end{align}
	Then, $x\mapsto V(t, x)$ is Lipschitz continuous for all $t\in [0, T]$, which, alongside the Lipschitz continuity of $t\mapsto V(t, x)$ for all $x\in \R$, allows us to conclude that $V$ is Lipschitz continuous on $[0, T]\times \R$.

	\ref{pr:V_C^12onC}\ \ 
	Since $V$ is continuous on $\cC$ (see \ref{pr:V_Lcontinuous}), $\mu$ in \eqref{eq:SDE} is Lipschitz continuous (actually, it suffices to require local $\alpha$-H\"older continuity) on $[0, T]\times\R$, and the diffusion coefficient is constant, one can borrow a classic result from the theory of parabolic partial differential equations \citep[see][Section 3, Theorem 9]{Friedman-1964-partial} to guarantee that, for an open rectangle $\cR\subset \cC$, the first initial-boundary value problem
	\begin{align}
	\mathbb{L}f - \lambda f &= 0 \qquad\text{in } \cR, \label{eq:PDE1} \\
	f &= V \qquad\! \text{on } \partial \cR, \label{eq:PDE2}
	\end{align}
	has a unique solution $f\in C^{1, 2}(\cR)$. Fix $(t, x)\in\cR$ and let $\tau_{\cR^c}$ be the first time $X_s^{t, x}$ exits $\cR$. Therefore, we can use the Itô formula on $f(X_s^{t,x})$ at $s = \tau_{\cR^c}$,
	together with \eqref{eq:PDE1} and \eqref{eq:PDE2}, to get the equality $\E_{t, x}[V(X_{t + \tau_{\cR^c}})] = f(t, x)$. Finally, notice that, due to the strong Markov property, $\E_{t, x}[V(X_{t + \tau_{\cR^c}})] = V(t, x)$, meaning that $f = V$ \citep[see][Section 7.1, for further details]{Peskir-2006-optimal}. The fact that $V$ is $C^{1, 2}$ on $\cD$ follows straightforwardly from $V = G$ on $\cD$.
	
	\ref{pr:x->V}\ \ We already mentioned that $x\mapsto V(t, x)$ is decreasing since $G$ is decreasing as well while~$x\mapsto X_s^{t, x}$ is increasing. For $(t, x)\in [0, T)\times \R$ such that $x\neq b(t)$, \eqref{eq:V_x} follows after recalling that $V$ is differentiable with respect to $x$ in $\cC$ and $\cD$, dividing by $\delta$ in \eqref{eq:Vx+delta-Vx}, and taking $\delta \rightarrow 0$, while using the dominated convergence theorem. The same procedure used in \eqref{eq:Vt-Vt-delta<} and \eqref{eq:Vt+delta-Vt>_1} yields \eqref{eq:V_t<} and \eqref{eq:V_t>}.
\end{proof}

Next, we look for smoothness of the free boundary in order to prove the smooth-fit condition, on which relies the uniqueness of the solution of the free-boundary problem associated to the OSP \eqref{eq:OSP}. The works of \cite{DeAngelis-2015-note}, \cite{Peskir-2019-continuity}, and \cite{DeAngelis-2019-Lipschitz} are a good compendium on the smoothness of the OSB. For time-homogeneous processes and smooth gain functions, \cite{DeAngelis-2015-note} provides the continuity of the free boundary for one-dimensional processes with locally Lipschitz-continuous drift and volatility. The two-dimensional case (including time-space diffusions) is addressed by \cite{Peskir-2019-continuity} in fairly general settings, which proves the impossibility of first-type discontinuities of the OSB for Mayer--Lagrange OSPs. 
Unfortunately,
continuity of the OSB is not enough to prove the smooth-fit condition.
Local Lipschitz continuity, however, suffices to derive the smooth-fit condition in our settings, as we prove later in Proposition~\ref{pr:smooth.fit}. 

In \cite{DeAngelis-2019-Lipschitz}, the local Lipschitz continuity of the free boundary in a high-dimensional framework is obtained, also for Mayer--Lagrange OSPs. However, some restrictive conditions are imposed that are not met in our setting, like $C^3$ smoothness of the gain function and a particular relation between the partial derivatives of $\InfGen G$ (see conditions (D), (F), and (G) in \cite{DeAngelis-2019-Lipschitz}).

The next proposition proves the local Lipschitz continuity of the OSB by adapting the methodology used by \cite{DeAngelis-2019-Lipschitz} to deal with our settings.

\begin{proposition}[Lipschitz continuity of $b$]\label{pr:OSB_LC}\ \\
	The OSB $b:[0,T]\to\mathbb{R}$ is Lipschitz continuous on any closed interval in $[0, T)$.
\end{proposition}

\begin{proof}
	Consider the function $W(t, x) = V(t, x) - G(x)$ and the closed interval $I = [\,\underline{t}, \bar{t}\,] \subset [0, T)$. Proposition \ref{pr:boundary_existence} guarantees that $b$ is bounded from above and, hence, we can choose $r > \sup \lrc{b(t): t\in I}$. 
	Hence, $I\times\lrc{r}\subset\cC$, and then $W(t, r) > 0$ for all $t\in I$. Since $W$ is continuous on $\overline{\cC}$ (see Proposition \ref{pr:V_regularity}), there exists a constant $a > 0$ such that $W(t, r) \geq a$ for all $t\in I$. Therefore, for all $\delta$ such that $0 < \delta \leq a$, the equation $W(t, x) = \delta$ has a solution in $\cC$ for all $t\in I$. This solution is unique for each $t$, as $\partial_x W > 0$ in $\cC$ (see \eqref{eq:V_x}). Hence we can denote it by $b_\delta(t)$, where $b_\delta:I\rightarrow (b(t), r]$. Away from the boundary, $W$ is regular enough to apply the implicit function theorem, which guarantees that $b_\delta$ is differentiable and
	\begin{align}\label{eq:b_delta'}
	b_\delta'(t) = -\partial_t W(t, b_\delta(t)) / \partial_x W(t, b_\delta(t)).
	\end{align}
	Notice that $b_\delta$ is increasing in $\delta$ and therefore it converges pointwise to some limit function $b_0$, which satisfies $b_0 \geq b$ in $I$ as $b_\delta > b$ for all $\delta$. Since $W(t, b_\delta(t)) = \delta$ and $W$ is continuous, it follows that $W(t, b_0(t)) = 0$ after taking $\delta\downarrow 0$, which means that $b_0 \leq b$ in $I$ and hence $b_0 = b$ in $I$.
	
	For $(t, x)\in \cC$ such that $t\in I$ and $x < r$, consider the stopping times $\tau^* = \tau^*(t, x)$ and
	$$
	\tau_r = \tau_r(t, x) := \inf\{s\geq 0 : X_{s}^{t, x} \notin I\times (-\infty, r) \}.
	$$ 
	Using \eqref{eq:partial_t_X<_1} from Lemma \ref{lm:bounds}, alongside the fact that $\tau^* \leq T - \underline{t}$, it readily follows after \eqref{eq:V_t<} that
	\begin{align}
		\partial_t W(t, x) \leq \Esp{\big|\partial_t X_{\tau^*}^{t, x}\big|} \leq L\lrp{1 + (T - \underline{t})e^{L(T - \underline{t})}}m(t, x), \label{eq:W_t<_1}
	\end{align}
	where $L$ is a positive constant coming from the Lipschitz continuity of $\partial_t\mu$ and $m$ is the function defined as
	\begin{align*}
	    m(t, x) := \Esp{\int_{0}^{\tau^*}\lrp{\left|X_{u}^{t, x}\right| + 1}\,\rmd u}.    
	\end{align*}
	Due to the tower property of conditional expectation and the strong Markov property, we have that
	\begin{align}
		m(t, x) 
		&= \Esp{\int_0^{\tau^*\wedge \tau_r}\lrp{\left|X_u^{t, x}\right| + 1}\,\rmd u + \Ind(\tau_r \leq \tau^*)m\left(t + \tau_r, X_{t + \tau_r}\right)}. \label{eq:m_split}
	\end{align}

	Notice that, for $b(t) < x < r$,  
	$\big(t + \tau_r, X_{t + \tau_r}^{t, x}\big)  \in \Gamma_t 
	:= \lrp{(t, \bar{t})\times\lrc{r}} \cup \lrp{\lrc{\bar{t}}\times(b(\bar{t}), r]}$ on the set $\lrc{\tau_r < \tau^*}$. Then, we get the following upper bound on $\lrc{\tau_r < \tau^*}$:
	\begin{align}
		m\left(t + \tau_r, X_{t + \tau_r}\right) \leq \sup_{(s, y) \in \Gamma_t}m\left(s, y\right) &= \sup_{(s, y) \in \Gamma_t} \Esp{\int_{0}^{\tau^*(s, y)}\lrp{\left|X_u^{s, y}\right| + 1}\,\rmd u} \nonumber \\
		&\leq T\sup_{(s, y) \in \Gamma_t} \sqrt{\E\bigg[\sup_{u\leq T - s}\lrp{\left|X_u^{s, y}\right| + 1}^2\bigg]} \nonumber \\
		&\leq T\sup_{(s, y) \in \Gamma_t} \sqrt{M_{T - s}^{(1)}\lrp{|y|^2 + 1}} \nonumber \\
		&\leq T\sqrt{M_{T}^{(1)}\lrp{\max\lrc{|b(\bar{t})|^2,|r|^2} + 1}}, \label{eq:m_bound}
	\end{align}
	where we used the Cauchy--Schwarz inequality and \eqref{eq:supX<} from Lemma \ref{lm:bounds}. Plugging \eqref{eq:m_bound} into \eqref{eq:m_split}, recalling \eqref{eq:W_t<_1}, and noticing that, for $u \leq \tau^*\wedge\tau_r$, $|X_u^{t, x}| \leq \max\{|\underline{b}|, r\}$ with $\underline{b} := \inf\{b(t) : t\in[0, T]\}$, we get that
	\begin{align}
		\partial_t W(t, x) \leq  K_I^{(1)}\Esp{\tau_r\wedge\tau^* + \Ind(\tau_r \leq \tau^*)} \label{eq:W_t<_2}
	\end{align}
	for some positive constant $K_I^{(1)}$.
	
	Using \eqref{eq:V_x} and \eqref{eq:1-partial_xX>} from Lemma \ref{lm:bounds}, and proceeding in the same way as for \eqref{eq:m_split}, we get the following lower bound for $\partial_x W(t, x)$:
	\begin{align}
	    \partial_x W(t, x) \geq M_{T - \underline{t}}^{(2)}\Esp{\tau^*} = M_{T - \underline{t}}^{(2)}\Esp{(\tau_r\wedge\tau^*) + \Ind(\tau_r \leq \tau^*)\Esp{\tau^*\lrp{t + \tau_r, X_{\tau_r}^{t, x}}}}. \label{eq:W_x>}
	\end{align}
	
	Now, by combining equations \eqref{eq:b_delta'}, \eqref{eq:W_t<_2}, and \eqref{eq:W_x>}, we obtain the following inequalities for a constant $K_I^{(2)} > 0$, $x_\delta = b_\delta(t)$, $\tau_\delta = \tau^*(t, x_\delta)$, and $\tau_r = \tau_r(t, x_\delta)$: 
	\begin{align}
		b_\delta'(t) \geq& - K_I^{(2)}\frac{\Esp{\tau_\delta\wedge\tau_r + \Ind(\tau_r\leq \tau_\delta)}}{\E\blrb{\tau_\delta\wedge\tau_r + \Ind(\tau_r\leq \tau_\delta)\E\blrb{\tau^*\blrp{t + \tau_r, X_{\tau_r}^{t, x_\delta}}}}} \nonumber \\
		\geq& -K_I^{(2)}\lrp{1 + \frac{\Prob{\tau_r\leq \tau_\delta}}{\E\blrb{\tau_\delta\wedge\tau_r + \Ind(\tau_r\leq \tau_\delta)\E\blrb{\tau^*\blrp{t + \tau_r, X_{\tau_r}^{t, x_\delta}}}}}} \nonumber \\
		\geq& -K_I^{(2)}\lrp{1 + \frac{\Prob{\tau_r\leq \tau_\delta, \tau_r = \bar{t} - t}}{\Esp{\tau_\delta\wedge\tau_r}} + \frac{\Prob{\tau_r\leq \tau_\delta, \tau_r < \bar{t} - t}}{\E\blrb{\Ind(\tau_r\leq \tau_\delta)\E\blrb{\tau^*\blrp{t + \tau_r, X_{\tau_r}^{t, x_\delta}}}}}} \nonumber \\
		\geq& -K_I^{(2)}\left(1 + \frac{\Prob{\tau_r\leq \tau_\delta, \tau_r = \bar{t} - t}}{\Esp{\Ind(\tau_r\leq \tau_\delta, \tau_r = \bar{t} - t)(\tau_\delta\wedge\tau_r)}}\phantom{\rule{0cm}{0.65cm}}\right.\nonumber\\ 
		&+ \left.\frac{\Prob{\tau_r\leq \tau_\delta, \tau_r < \bar{t} - t}}{\E\blrb{\Ind(\tau_r\leq \tau_\delta, \tau_r < \bar{t} - t)\E\blrb{\tau^*\blrp{t + \tau_r, X_{\tau_r}^{t, x_\delta}}}}}\right) \nonumber \\
        \geq& -K_I^{(2)}\lrp{1 + \frac{1}{\bar{t} - t} + \frac{1}{\inf_{t\in I}\Esp{\tau^*(t, r)}}}, \label{eq:b'>}
	\end{align}
	where in the last inequality we used that $X_{\tau_r}^{t, x_\delta} = r$ in the set $\lrc{\tau_r\leq \tau_\delta, \tau_r < \bar{t} - t}$.
	
	Since $I\times\lrc{r}\subset\cC$, there exists $\varepsilon > 0$ such that $\cR_\varepsilon := [\,\underline{t}, \bar{t} + \varepsilon]\times(r - \varepsilon, r + \varepsilon)\subset\cC$. Notice that $\tau^*(t, r) > \tau_\varepsilon(t, r)$ for all $t\in I$, with $\tau_\varepsilon(t, r) = \inf\blrc{u\geq 0 : \blrp{t + u, X_{u}^{t, r}} \notin \cR_\varepsilon}$. Hence,
	\begin{align*}
	    \inf_{t\in I} \Esp{\tau^*(t, r)} &\geq \inf_{t\in I}\Esp{\tau_\varepsilon(t, r)} \\
	    &\geq \inf_{t\in I} \bigg\{(\bar{t} + \varepsilon - t)\mathbb{P}\bigg(\sup_{u\leq \bar{t} + \varepsilon - t} \lrav{X_u^{t, r} - r} < \varepsilon\bigg)\bigg\} \\
	    &\geq \varepsilon \inf_{t\in I} \mathbb{P}\bigg(\sup_{u\leq \bar{t} + \varepsilon - t} \lrav{X_u^{t, r} - r} < \varepsilon\bigg) =: K_{I,\varepsilon} > 0.
	\end{align*}
	
	Therefore, going back to \eqref{eq:b'>}, we have that, for all $t\in I' = [\,\underline{t}, \bar{t} - \varepsilon]$ and $\varepsilon > 0$ small enough,
	\begin{align}
		b_\delta'(t) \geq -K_I^{(2)}\lrp{1 + \varepsilon^{-1} + K_{I,\varepsilon}}. \label{eq:b'>const_2}
	\end{align}
	
	To find a bound in the opposite direction that is also uniform with respect to $\delta$ and for all $t\in I'$, we consider \eqref{eq:V_t>} and \eqref{eq:W_t<_1} and use the Markov inequality to get
	\begin{align}
		\partial_t W(t, x_\delta) \geq& -L\lrp{1 + (T - \underline{t})e^{L(T - \underline{t})}}\Esp{\int_{0}^{\tau_\delta}\lrp{\left|X_u^{t, x_\delta}\right| + 1}\,\rmd u}\nonumber\\
		&- \widetilde{L}(T - t)^{-1}\Esp{\tau_\delta} \label{eq:W_t>}
	\end{align}
	for some positive constant $\widetilde{L}$. Hence, relying on the same arguments as those used to get \eqref{eq:b'>const_2}, but with \eqref{eq:W_t>} instead of \eqref{eq:W_t<_2}, we obtain  
	\begin{align}
		b_\delta'(t) \leq K\lrp{1 + \varepsilon^{-1} + K_{I,\varepsilon}} + \widetilde{L}(T - \bar{t})^{-1}. \label{eq:b'<const}
	\end{align}
	
	We have proved that $|b_\delta'(t)|$ is bounded by a constant, uniformly in $\delta$ and for all $t\in I'$. Then the Arzelà--Ascoly theorem gives that $b_\delta$ converges uniformly to $b$ in $I'$ as $\delta\rightarrow 0$, which implies that $b$ is Lipschitz continuous on $I'$.
\end{proof} 

The next proposition gives the smooth-fit condition for our OSP. Its proof relies on Lemma~\ref{lm:LIL} and the fact that $b$ is locally Lipschitz continuous. These two combined results imply that the process $X^{t, x}$ enters the interior of $\cD$ immediately for $x = b(t)$ (see Remark 4.5 in \cite{DeAngelis-2019-Lipschitz}). The work in \cite{DeAngelis-2020-global}, specifically Corollary 6, is then used to obtain $\tau^*(t, b(t)+\varepsilon) \rightarrow 0$ $\Pr$-a.s. when $\varepsilon\downarrow 0$, from which the smooth-fit condition follows.

\begin{proposition}[Smooth-fit condition]\label{pr:smooth.fit}\ \\
	The smooth-fit condition holds, i.e., $\partial_xV(t,b(t)) = -1$ for all $t\in[0,T)$. 
\end{proposition}

\begin{proof}
	Take a point $(t, b(t))$ for $t\in[0, T)$ and consider $\delta > 0$. Since $(t, b(t)) \in \cD$ and $(t, b(t) + \delta)~\in~\cC$, we get that $\delta^{-1}(V(t, b(t) + \delta) - V(t, b(t))) \geq \delta^{-1}(G(b(t) + \delta) - G(b(t))) = -1$. Therefore, $\partial_x^+V(t, b(t)) \geq -1$. 
	
	Besides, reasoning as in \eqref{eq:Vx+delta-Vx}, we get that
	\begin{align*}
	\delta^{-1}(V(t, b(t) + \delta) - V(t, b(t))) &\leq -\E\left[e^{-\lambda\tau_\delta}\partial_xX_{\tau_\delta}^{t, x_\delta}\right] \\
	&\leq \sup_{x\in (b(t), A)}\E\bigg[\sup_{s\leq T - t}\partial_x X_s^{t, x}\bigg] < \infty,
	\end{align*}
	where $\tau_\delta = \tau^*(t, b(t) + \delta)$ and $x_\delta \in(b(t), b(t) + \delta)$. From \eqref{eq:SDE_x^}, it is easy to see that the supremum is finite (actually, it is lower than~$1$) since $\partial_x\mu < 0$. Then, we can apply the dominated convergence theorem to obtain
	\begin{align}\label{eq:Vx-smooth.fit<}
            \partial_x^+V(t, b(t)) \leq -\E\left[e^{-\lambda\tau_0}\partial_xX_{\tau_0}^{t, b(t)}\right]
	\end{align}
	with $\tau_0 := \lim_{\delta\rightarrow0} \tau_\delta$ ($\tau_0$ is well-defined since the sequence ${\tau_\delta}$ decreases with respect to $\delta$ $\Pr$-a.s.).

    We now prove that $\Pro{\tau_0 = 0}{t, b(t)} = 1$. We do so by summoning Corollary 6 from \cite{DeAngelis-2020-global} and proving that $(t, b(t))$ is probabilistically regular for $\interior{\cD}$, that is, $\Pro{\tau_{\interior{\cD}} = 0}{t, b(t)} = 1$, where
    \begin{align*}
        \tau_{\interior{\cD}} := \inf\lrc{u \geq 0 : (t+u, X_{t+u}) \in \interior{\cD}}.
    \end{align*}
    Indeed,
    \begin{align*}
        \Pro{\tau_{\interior{\cD}} = 0}{t, b(t)} &= \lim_{\varepsilon\downarrow 0} \Pro{\tau_{\interior{\cD}} < \varepsilon}{t, b(t)} \\
        &= \lim_{\varepsilon\downarrow 0}\Pro{\inf_{u\in(0,\varepsilon)} \lrp{X_{t+u} - b(t+u)} < 0}{t, b(t)} \\
        &= \lim_{\varepsilon\downarrow 0}\Pro{\inf_{u\in(0,\varepsilon)} \frac{X_{t+u} - b(t+u)}{\sqrt{2u\ln(\ln(1/u))}} < 0}{t, b(t)} \\
        &\geq \lim_{\varepsilon\downarrow 0}\Pro{\inf_{u\in(0,\varepsilon)} \frac{X_{t+u} - b(t) + Lu}{\sqrt{2u\ln(\ln(1/u))}} < 0}{t, b(t)} \\
        &= \Pro{\liminf_{u\downarrow 0} \frac{X_{t+u} - b(t) + Lu}{\sqrt{2u\ln(\ln(1/u))}} < 0}{t, b(t)} = 1,
    \end{align*}
    where, in the first inequality, $L$ is a positive constant that comes from the local Lipschitz continuity of $b$ (see Proposition \ref{pr:OSB_LC}), and the last equality holds due to \eqref{eq:OU-LIL_inf}. Hence, \eqref{eq:Vx-smooth.fit<} yields $\partial_x^+V(t, b(t)) \leq -1$ and, consequently, $\partial_x^+V(t, b(t)) = -1$.
	
	The smooth-fit condition arises by recalling that $\partial_x^-V(t, b(t)) = -1$ since $V = G$ on~$\cD$.   
\end{proof}

We conclude this section with two technical lemmas. One on the boundedness of several functionals of the underlying process $X$ that have been used throughout the previous proofs, and another on a law of the iterated logarithm for $X$. 

\begin{lemma}[Some bounds on functionals of $X$] \label{lm:bounds}\ \\
	The following inequalities hold for positive constants $L$, $M_s^{(i)}$, $i = 1, 2, 3$, and $(t, x) \in [0, T)\times \R$ and $s\in(0, T - t]$:
	\begin{align}
		\Esp{\sup_{u\leq s}\lrp{\left|X_u^{t, x}\right| + 1}^2} &\leq M_s^{(1)}\lrp{|x|^2 + 1}, \label{eq:supX<} \\
		\partial_x X_{s}^{t, x}  &\leq 1 - M_s^{(2)}s, \label{eq:1-partial_xX>}\\
		\Esp{\sup_{u\leq s}\left|\partial_tX_u^{t, x}\right|} &\leq L\lrp{1 + se^{Ls}}\Esp{\int_{0}^{s}\lrp{\left|X_u^{t, x}\right| + 1}\,\rmd u}  \label{eq:partial_t_X<_1} \\
		&\leq M_s^{(3)}\lrp{|x| + 1}. \label{eq:partial_t_X<}
	\end{align}
	Moreover, $s\mapsto M_s^{(i)}$ is an increasing function for $i = 1, 3$, and is decreasing for $i = 2$.
\end{lemma}

\begin{proof}
	Since $(a + b + c)^2 \leq 3(a^2 + b^2 + c^2)$ and due to the Lipschitz continuity of $x\mapsto\mu(t, x)$,
	\begin{align*}
	\lrp{\left|X_s^{t, x}\right| + 1}^2 &\leq 3\lrp{\left|x\right| + 1}^2 + 3\int_0^s \lrp{\mu(t + r, X_r^{t, x})}^2\,\rmd r + 3 \lrp{W_s}^2 \\
	&\leq 3\lrp{\left|x\right| + 1}^2 + 3L\int_0^s \lrp{\left|X_r^{t, x}\right| + 1}^2\,\rmd r + 3\lrp{W_s}^2
	\end{align*}
	for some positive constant $L$. Hence, using the maximal inequalities of Theorem 14.13(d) from \cite{Schilling-2012-Brownian}, it follows that
	\begin{align*}
	\Esp{\sup_{u\leq s}\lrp{\left|X_u^{t, x}\right| + 1}^2} &\leq 3\lrp{\left|x\right| + 1}^2 + 3L\int_0^s \Esp{\lrp{\left|X_r^{t, x}\right| + 1}^2}\,\rmd r + 3\Esp{\sup_{u\leq s}\lrp{W_u}^2} \\
	&\leq 3\lrp{\left|x\right| + 1}^2 + 12 s + 3L\int_0^s \Esp{\sup_{u\leq r}\lrp{\left|X_u^{t, x}\right| + 1}^2}\,\rmd r.
	\end{align*}
	Therefore, Gronwall's inequality \cite[Theorem A.43]{Schilling-2012-Brownian} guarantees that
	\begin{align*}
	\Esp{\sup_{u\leq s}\lrp{\left|X_u^{t, x}\right| + 1}^2} &\leq 3\lrp{\left|x\right| + 1}^2 + 12 s + 3L\int_{0}^{s}\lrp{3\lrp{\left|x\right| + 1}^2 + 12 u}e^{3L(s - u)}\,\rmd u \\
	&\leq 3\lrp{\left|x\right| + 1}^2 + 12 s + 3Le^{3Ls}\lrp{3\lrp{\left|x\right| + 1}^2s + 6 s^2}, 
	\end{align*}
	from which \eqref{eq:supX<} follows.
	
	To obtain \eqref{eq:1-partial_xX>}, we use that $\mu$ is Lipschitz continuous and $\partial_x\mu < 0$, together with representation \eqref{eq:SDE_x^} and \eqref{eq:V_x}, leading to
	\begin{align*}
		1 - \partial_x X_{s}^{t, x} &= \int_{0}^{s}-\partial_x\mu(t + u, X_u^{t, x})\partial_x X_u^{t, x}\,\rmd u  \\
		&= \int_{0}^{s}-\partial_x\mu(t + u, X_u^{t, x})\exp\lrc{\int_{0}^{u}\partial_x\mu(t + r, X_r^{t, x})\,\rmd r}\,\rmd u  \\
		&\geq \exp\lrc{\int_{0}^{s}\partial_x\mu(t + u, X_u^{t, x})\,\rmd u}\int_{0}^{s}-\partial_x\mu(t + u, X_u^{t, x})\,\rmd u  \\
		&\geq e^{-Ls}\int_{0}^{s}-\partial_x\mu(t + u, X_u^{t, x})\,\rmd u  \\
		&\geq e^{-Ls}\min_{u\in[0, T]}\lrc{-\bar{\mu}(u)}s,
	\end{align*}
	where $L$ is the Lipschitz constant. 
	
	Let us prove now \eqref{eq:partial_t_X<}. To do so, notice first that \eqref{eq:SDE_t} implies
	\begin{align*}
		\left|\partial_t X_s^{t, x}\right| &\leq a_1(s, t, x) + \int_{0}^{s}a_2(u, t, x)\left|\partial_t X_u^{t, x}\right|\,\rmd u,
	\end{align*}
	with 
	\begin{align*}
		a_1(s, t, x) := \int_{0}^{s}\left|\partial_t\mu(t + u, X_u^{t, x})\right|\,\rmd u, \quad
		a_2(u, t, x) := \left|\partial_x\mu(t + u, X_u^{t, x})\right|.
	\end{align*}
	Therefore, an application of Gronwall's inequality yields
	\begin{align*}
		\left|\partial_t X_s^{t, x}\right| &\leq a_1(s, t, x) + \int_{0}^{s}a_1(u, t, x)a_2(u, t, x)\exp\left\{\int_u^s a_2(r, t, x)\,\rmd r\right\}\,\rmd u.
	\end{align*}
	Hence, due to the Lipschitz continuity of $\mu$ and $\partial_t\mu$, and using \eqref{eq:supX<} alongside the Cauchy--Schwarz inequality, we have that
	\begin{align}
		\Esp{\sup_{u\leq s}\left|\partial_t X_{u}^{t, x}\right|}  &\leq \Esp{a_1(s, t, x)\lrp{1 + \int_{0}^{s}a_2(u, t, x)\,\rmd \nonumber u\exp\left\{\int_{0}^{s}a_2(u, t, x)\,\rmd u\right\}}} \nonumber \\
		&\leq L\lrp{1 + se^{Ls}}\Esp{\int_{0}^{s}\lrp{\left|X_u^{t, x}\right| + 1}\,\rmd u} \nonumber \\
		&\leq L\lrp{1 + se^{Ls}}s\sqrt{\Esp{\sup_{u\leq s}\lrp{\left|X_u^{t, x}\right| + 1}^2}} \nonumber \\
		&\leq L\lrp{1 + se^{Ls}}s\sqrt{M_s^{(1)}\lrp{|x|^2 + 1}}, \nonumber
	\end{align}
	which entails \eqref{eq:partial_t_X<_1} and \eqref{eq:partial_t_X<}.
\end{proof}

\begin{lemma}[Law of the iterated logarithm for $X$]\label{lm:LIL}
    The following relations hold for all $(t, x)\in[0, T]\times\R$:
    \begin{align}\label{eq:OU-LIL_sup}
        \Pro{\limsup_{u\downarrow0}\frac{X_{t+u} - x}{\sqrt{2u\ln(\ln(1/u))}} = 1}{t, x} = 1,
    \end{align}
    and
    \begin{align}\label{eq:OU-LIL_inf}
        \Pro{\liminf_{u\downarrow0}\frac{X_{t+u} - x}{\sqrt{2u\ln(\ln(1/u))}} = -1}{t, x} = 1.
    \end{align}
\end{lemma}

\begin{proof}
    Let $m_{t, x}(u) := \Es{X_{t+u}}{t, x}$, $f_t(u) = \exp\big\{-\int_t^{t+u}\theta(r)\,\rmd r\big\}$, and $h_t(u) = \int_t^{t+u} f_t^{-2}(u) \,\rmd u$. It is known (see, e.g., \cite{Buonocore-2013-some}) that, under $\Pr_{t, x}$, $X_{t+u}$ admits the representation $X_{t+u} = m_{t, x}(u) + f_t(u)B_{h_t(u)}$, where $(B_t)_{t\in\R_+}$ is a standard Brownian motion. Hence, the law of the iterated logarithm for the Brownian motion establishes that
    \begin{align*}
        \Pro{\limsup_{u\rightarrow0}\frac{(X_{t+u} - m_{t, x}(u))/f_t(u)}{\sqrt{2h_t(u)\ln(\ln(1/h_t(u)))}} = 1}{t, x} = 1.
    \end{align*}
    Then, \eqref{eq:OU-LIL_sup} is proved after realizing that $m_{t, x}(u)\rightarrow x$, $f_t(u)\rightarrow 1$, and $h_t(u)/u \rightarrow 1$, as $u\downarrow0$. Relation \eqref{eq:OU-LIL_inf} follows identical steps.
\end{proof}

\section{The value formula and the free-boundary equation}\label{sec:pricing_formula}

Our main result is collected in the following theorem. It gives a characterization of the OSB as the unique solution, up to regularity conditions, of a type-two Volterra integral equation. It also provides a formula for the value function, which can be regarded as the fair price of the American option under the natural measure of the OU process. This contract might find application from the point of view of small investors who, due to their limited capital and impossibility to trade in continuous time (also due to transaction costs), cannot perform hedging strategies to reproduce the conditions of a risk-free measure settings under which the valuation of options is typically done to keep the market free of arbitrage opportunities. 

The value formula is given in terms of the early-exercise premium representation. That is, the sum of the price of the European put option written on the same asset and expiring on the same date, and the early-exercise premium, i.e., the cost of being able to exercise the option before maturity.

\begin{theorem}[Free-boundary equation and value formula]\label{th:fb}\ \\
	The OSB related to the OSP \eqref{eq:OSP} satisfies the integral equation
	\begin{align}
        b(t) =&\; A - K_\lambda(A, 1, t, b(t), T, A) \nonumber\\
        &- \int_t^T K_\lambda(\lambda A + \theta(u)\alpha(u), \lambda + \theta(u), t, b(t), u, b(u))\,\rmd u,\label{eq:free-boundary_equation_gaussian}
    \end{align}
	where, for $c_1, c_2, x_1, x_2 \in\R$ and $t_1, t_2 \in[0, T]$ such that $t_2 \geq t_1$,
	\begin{align}
        K_\lambda(c_1, c_2, t_1, x_1, t_2, x_2) :=&\; e^{-\lambda(t_2 - t_1)}\Es{\lrp{c_1 - c_2 X_{t_2}}\Ind(X_{t_2}\leq x_2)}{t_1, x_1} \label{eq:truncated_normal_mean} \\
        =&\; e^{-\lambda(t_2 - t_1)}\bigg\{(c_1 - c_2\nu(t_1, x_1, t_2))\Phi\lrp{\frac{x_2 - \nu(t_1, x_1, t_2)}{\gamma(t_1, t_2)}} \nonumber\\
        &+c_2\gamma(t_1, t_2)\phi\lrp{\frac{x_2 - \nu(t_1, x_1, t_2)}{\gamma(t_1, t_2)}}\bigg\} \nonumber,
    \end{align}
	where $\phi$ and $\Phi$ are the density and distribution functions of a standard normal, and $\nu$ and $\gamma$ depend on the time-dependent OU parameters:
	\begin{align}
        \nu(t_1, x, t_2)&= \exp\lrc{-\int_{t_1}^{t_2}\theta(r)\,\rmd r}x + \int_{t_1}^{t_2} \exp\lrc{-\int_{r}^{t_2}\theta(s)\,\rmd s}\theta(r)\alpha(r)\,\rmd r, \label{eq:mean} \\ 
        \gamma^2(t_1, t_2) &= \int_{t_1}^{t_2} \exp\lrc{-2\int_r^{t_2}\theta(s)\,\rmd s}\,\rmd r. \label{eq:variance}
    \end{align}
    Moreover, \eqref{eq:free-boundary_equation_gaussian} has a unique solution among the class of continuous functions of bounded variation $c:[0,T]\rightarrow\mathbb{R}$ that satisfy $c(t)< A$ for all $t\in[0, T)$.
    
    Also, the associated value function is given by
    \begin{align}\label{eq:pricing_formula_gaussian}
        V\left(t, x\right) =&\; K_\lambda(A, 1, t, x, T, A) \\
        &+ \int_t^T K_\lambda(\lambda A + \theta(u)\alpha(u), \lambda + \theta(u), t, x, u, b(u))\,\rmd u. \nonumber
    \end{align}
\end{theorem}

\begin{proof}
    Propositions \ref{pr:boundary_existence}--\ref{pr:smooth.fit} allow applying an extension of the Itô formula to $e^{-\lambda s}V(t + s, X_{t + s})$. This extension is derived in \cite{Peskir-2005-change} and is reformulated in Lemma A2 from \cite{DAuria-2020-discounted} in a way that fits our settings more straightforwardly. After setting $s = T - t$, taking $\Pr_{t, x}$-expectation, and shifting the integrating variable $t$ units, it follows that
    \begin{align}\label{eq:V_after_Ito}
    	V(t, x) = e^{-\lambda (T - t)}\E_{t, x}\left[G(X_T)\right] - \int_{t}^{T}\E_{t, x}\left[e^{-\lambda (u - t)}(\mathbb{L} V - \lambda V)(u, X_u)\right]\,\rmd u,
    \end{align} 
    where the martingale term vanishes after taking $\Pr_{t, x}$-expectation, and the local time term does not appear due to the smooth-fit condition. Recall that $\mathbb{L} V = \lambda V$ on $\cC$ and $V = G$ on $\cD$. Therefore, \eqref{eq:V_after_Ito} turns into the value formula
    \begin{align}\label{eq:pricing_formula_general}
    	V\left(t, x\right) =&\; e^{-\lambda (T - t)}\E_{t, x}\left[\left(A - X_T\right)\Ind\left(X_T \leq A\right)\right] \\
    	&+ \int_{t}^{T}\E_{t, x}\left[e^{-\lambda (u - t)}\left(\lambda\left(A - X_{u}\right) + \mu(u,                 X_u)\right)\Ind\left(X_{u}\leq b(u)\right)\right]\,\rmd u. \nonumber 
    \end{align}
    
    By taking $x \uparrow b(t)$ in \eqref{eq:pricing_formula_general} we get the free-boundary equation
    \begin{align}\label{eq:free-boundary_equation_general}
	    b(t) =&\; A - e^{-\lambda (T - t)}\E_{t, b(t)}\left[\left(A - X_T\right)\Ind\left(X_T \leq A\right)\right] \\
	    & - \int_{t}^{T}\E_{t, b(t)}\left[e^{-\lambda (u - t)}\left(\lambda\left(A - X_{u}\right) + \mu(u, X_u)\right)\Ind\left(X_{u}\leq b(u)\right)\right]\,\rmd u. \nonumber
    \end{align}
    
    The uniqueness of the solution of \eqref{eq:free-boundary_equation_general}, up to the conditions stated in the theorem, follows well-known arguments derived in \cite[Theorem 3.1]{Peskir-2005-American}.

    We can provide a more tractable expression for both the value formula \eqref{eq:pricing_formula_general} and the free-boundary equation \eqref{eq:free-boundary_equation_general} by taking advantage of the Gaussianity of $X$. Indeed, it is easy to get \citep[see, e.g.,][Section 4.4.9]{Gardiner-2004-handbook} that, under $\Pr_{t, x}$, $X_u$ has a normal distribution with mean and variance given by \eqref{eq:mean} and \eqref{eq:variance} for all $u\in[t, T]$. This leads, after some algebraic manipulation, to equations \eqref{eq:free-boundary_equation_gaussian} and \eqref{eq:pricing_formula_gaussian}.
\end{proof}

Further flexibility of the underlying process can be achieved by allowing time-dependent volatility, like in the popular Hull--White model. We show in the next remark that such an extension requires no extra analysis, as it can be reduced to the constant-volatility model \eqref{eq:SDE} by means of a time change. 

\begin{remark}[Time-dependent deterministic volatility]\label{rmk:time_dependent_volatility}\ \\
    Let $X_t(\theta, \alpha, \sigma, T) = \lrp{X_t}_{t \in [0, T]}$ be a stochastic process satisfying the stochastic differential equation
\begin{align}\label{eq:GM_SDE}
    \rmd X_t(\theta, \alpha, \sigma, T) = \theta(t)(\alpha(t) - X_t(\theta, \alpha, \sigma, T))\,\rmd t + \sigma(t)\,\rmd W_t,\ \ 0 \leq t \leq T,
\end{align}
where $\theta, \alpha, \sigma:[0, T]\rightarrow \R$ are continuously differentiable functions such that $\theta(t) > 0$ and $\sigma(t) > 0$ for all $t\in[0, T]$, and $\int_0^T\sigma^2(t)\,\rmd t < \infty$. Define the time change $\tilde{t} = \tilde{h}(t) = \int_0^t \sigma^2(u)\,\rmd u$ and the inverse transformation $h = \tilde{h}^{-1}$. Consider the time-changed process $\wt{X} = \big(\wt{X}_{\tilde{t}}\big)_{\tilde{t}\in[0, \wt{T}]}$ defined by $\wt{X}_{\tilde{t}} = X_{h(\tilde{t})}(\theta, \alpha, \sigma, T)$ and $\wt{T} = \tilde{h}(T)$. Then,
\begin{align*}
    \wt{X}_{\tilde{t}} &= X_0(\theta, \alpha, \sigma, T) + \int_0^{h(\tilde{t})} \theta(u)(\alpha(u) - X_u(\theta, \alpha, \sigma, T))\,\rmd u + \int_0^{h(\tilde{t})}\sigma(u)\,\rmd W_u.
\end{align*}    
The process $\wt{W}_{\tilde{t}} = \int_0^{h(\tilde{t})}\sigma(u)\,\rmd W_u$ is a standard Brownian motion. Hence, $\wt{X}$ solves
\begin{align*}
    \rmd \wt{X}_{\tilde{t}} &= h'(\tilde{t})\theta(h(\tilde{t}))\big(\alpha(h(\tilde{t})) - X_{h(\tilde{t})}(\theta, \alpha, \sigma, T)\big)\,\rmd \tilde{t} + \rmd\wt{W}_{\tilde{t}} \\
    &= \wt{\theta}(\tilde{t})\big(\wt{\alpha}(\tilde{t}) - \wt{X}_{\tilde{t}}\big)\,\rmd \tilde{t} + \rmd\wt{W}_{\tilde{t}},
\end{align*}
where $\wt{\theta}(\tilde{t}) = \sigma^{-2}(h(\tilde{t}))\theta(h(\tilde{t}))$ and $\wt{\alpha}(\tilde{t}) = \alpha(h(\tilde{t}))$, 
that is, $\wt{X}$ is a time-dependent OU process with constant volatility. Hence, we can rely on equations \eqref{eq:pricing_formula_gaussian} and \eqref{eq:free-boundary_equation_gaussian} to compute the value function and OSB associated with $\wt{X}$, $\wt{V}$, and $\tilde{b}$, from where we can simply obtain the value function and OSB related to $X(\theta, \alpha, \sigma, T)$ as $V(t, x) = \wt{V}(h(\tilde{t}), x)$ and $b(t) = \tilde{b}(h(\tilde{t}))$.
\end{remark}

The optimal exercise of the put and the call American options is fundamentally the same problem due to their symmetry with respect to the strike price, as shown in the next remark.

\begin{remark}[Put-call parity]\label{rmk:put-call}\ \\
    Let $V_\mathrm{p}$ be the value function associated with an American put option written on the time-dependent OU processes $\smash{X^{(\mathrm{p})} = \blrp{X_t^{(\mathrm{p})}}_{t\in[0, T]}}$ and with strike price $A$. Then,
    \begin{align*}
        V_\mathrm{p}(t, x) :=&\; \sup_{\tau\leq T - t}\Es{e^{-\lambda\tau}\lrp{A - X_{t+\tau}^{(\mathrm{p})}}^+}{t, x} \\
        =&\; \sup_{\tau\leq T - t}\Es{e^{-\lambda\tau}\lrp{\lrp{2A - X_{t+\tau}^{(\mathrm{p})}} - A}^+}{t, x} \\
        =&\; \sup_{\tau\leq T - t}\Es{e^{-\lambda\tau}\lrp{X_{t+\tau}^{(\mathrm{c})} - A}^+}{t, 2A - x} =: V_\mathrm{c}(t, 2A - x), 
    \end{align*}
    where $\smash{X^{(\mathrm{c})} = \blrp{X_t^{(\mathrm{c})}}_{t\in[0, T]}}$ is the time-dependent OU process such that $X^{(\mathrm{c})} = 2A - X^{(\mathrm{p})}$, and $V_\mathrm{c}$ is the value function of the American call option with strike price $A$ and written on $X^{(\mathrm{c})}$. Similarly, the optimal exercise strategy for the call option is to stop the process $X^{(\mathrm{c})}$ the first time it is below the function $b_\mathrm{c}:[0, T]\rightarrow\R$ such that $b_\mathrm{c}(t) = 2A - b_\mathrm{p}(t)$, where $b_\mathrm{p}$ is the OSB associated to $V_\mathrm{p}$. 
\end{remark}

\section{Numerical experiments}\label{sec:numerical_results}

Non-linear Volterra-type integral equations are hard to solve, and the free-boundary equation \eqref{eq:free-boundary_equation_gaussian} is not an exception. Mainly, there are two methods to solve these types of integral equations arising from OSPs, both relying on the contraction principle: (1) a backward induction approach that recursively builds the boundary by leveraging that $b(t)$ depends only on $\lrc{b(s)}_{s\in[t, T]}$; and~(2) the well-known method of Picard iterations that, given an initial candidate boundary $b^{0}$, consequently computes subsequent boundaries until a convergence criterion is fulfilled. See \citet[Chapter 8]{Detemple-2005-American-style} for more details on the method of backward induction and \cite{Peskir-2002-nonlinear} for implementations. The Picard iteration scheme has been used by, e.g., \cite{Detemple-2020-value} and \cite{DeAngelis-2020-optimal}. As far as we know, the convergence of both approaches has not been formally addressed in general settings, and authors tend to provide numerical evidence of the error decreasing as the number of iterations grows. This happens because the kernel of the integral equation is typically highly non-linear and, hence, it becomes challenging to prove the contracting mapping principle for the integral operator that characterizes the OSB, which is the most obvious and immediate strategy to follow to get the convergence of Picard's iteration algorithm. We adopted the Picard method to display visual insight into the OSB's shape, as our numerical experiments suggested similar accuracy and faster computations compared to the backward induction approach (within the set of explored settings).

Take the partition $0 = t_0 < t_1 < \cdots < t_N = T$ for $N\in\mathbb{N}$. We initialize the Picard scheme by starting with the constant boundary $b^{(0)}(t) = b(T)$ for all $t\in[0, T]$. Recall from Proposition~\ref{pr:boundary_existence} that $b(T) = \min\{A, \gamma(T)\}$ is known, with $\gamma(T)$ as in \eqref{eq:gamma}. Each iteration of the algorithm updates the boundary according to the following formula, which is a right Riemann sum version of the integral in \eqref{eq:free-boundary_equation_gaussian}:
\begin{align}
    b^{(k)}(t_i) =&\; A - K_\lambda\blrp{A, 1, t_i, b^{(k - 1)}(t_i), T, A} \label{eq:Picard}\\
    & - \sum_{j = i + 1}^{N} K_\lambda\blrp{\lambda A + \theta(t_j)\alpha(t_j), \lambda + \theta(t_j), t_i, b^{(k - 1)}(t_i), t_j, b^{(k - 1)}(t_j)}. \nonumber
\end{align}
We stop the iterations of the fixed-point algorithm \eqref{eq:Picard} at the first $k = 1, 2, \ldots$ such that the (squared) distance $d_k := \sum_{i=0}^N \blrp{b^{(k)}(t_i) - b^{(k-1)}(t_i)}^2$ is less than $\delta=10^{-3}$. Computational experiments suggested that a mesh that is denser at the terminal time $T$ performs better and that it is preferable that the distances between consecutive nodes narrow smoothly. For that reason, we used the logarithmically-spaced partition $t_i = \ln\lrp{1 + i(e - 1)/N}$ with $N = 200$ (unless otherwise specified).

For Figures \ref{fig:OSB_shape}--\ref{fig:OSB_convergence_and_discount}, the images on top show computer drawings of the obtained OSB from the Picard iteration \eqref{eq:Picard}. The images below represent the sequences of errors of the algorithm, that is, the $x$-axis accounts for the iteration number $k = 1, 2, \ldots$ and the $y$-axis for $d_k$. The decrease in these error sequences empirically corroborates the convergence of the algorithm.

\begin{figure}
\centering
\subfloat[$\theta(t) = 1,\, \sigma(t) = 1$.]{%
\resizebox*{0.32\textwidth}{!}{\includegraphics{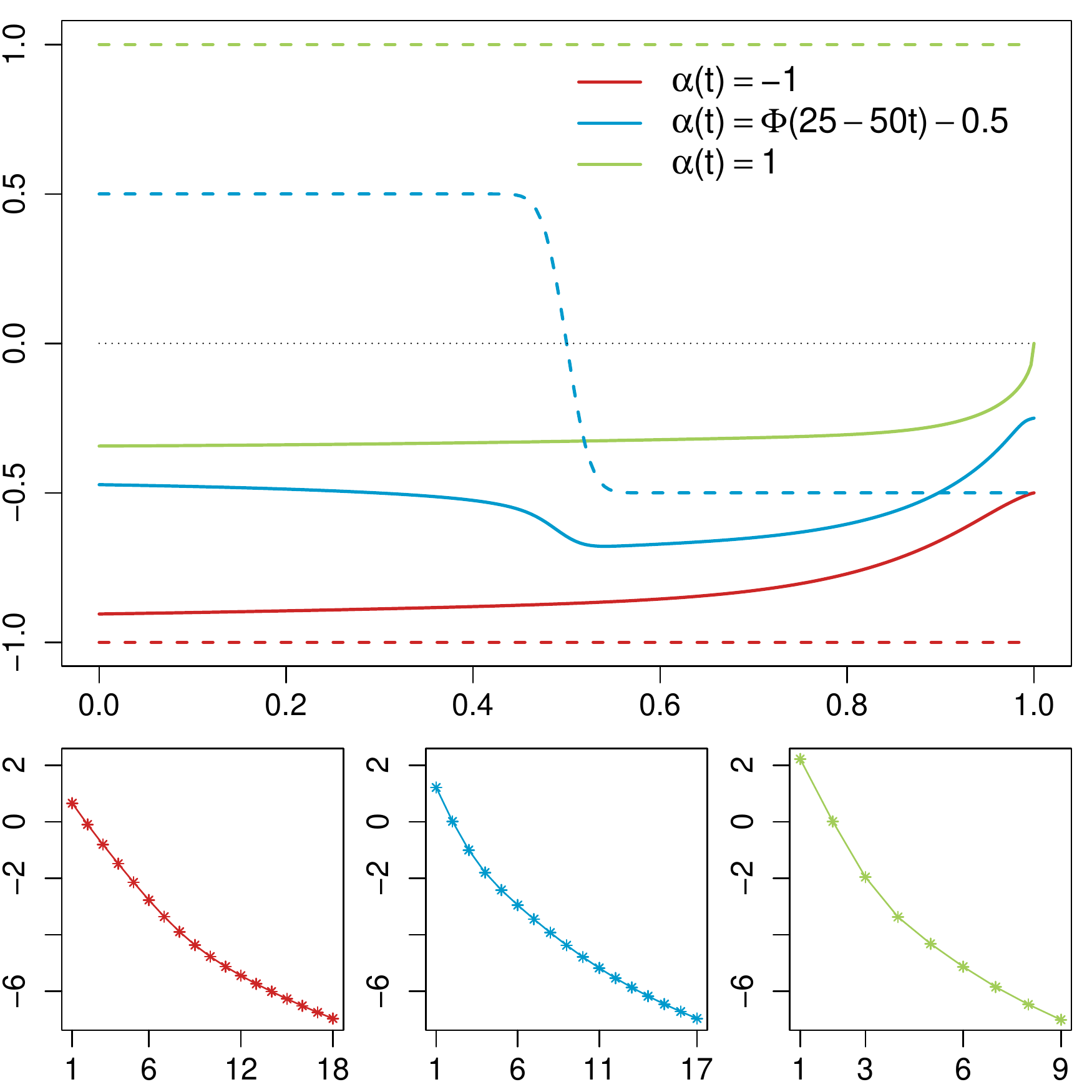}}}\hspace{2pt}
\subfloat[$\sigma(t) = 1,\, \alpha(t) = e^{0.5t}$.]{%
\resizebox*{0.32\textwidth}{!}{\includegraphics{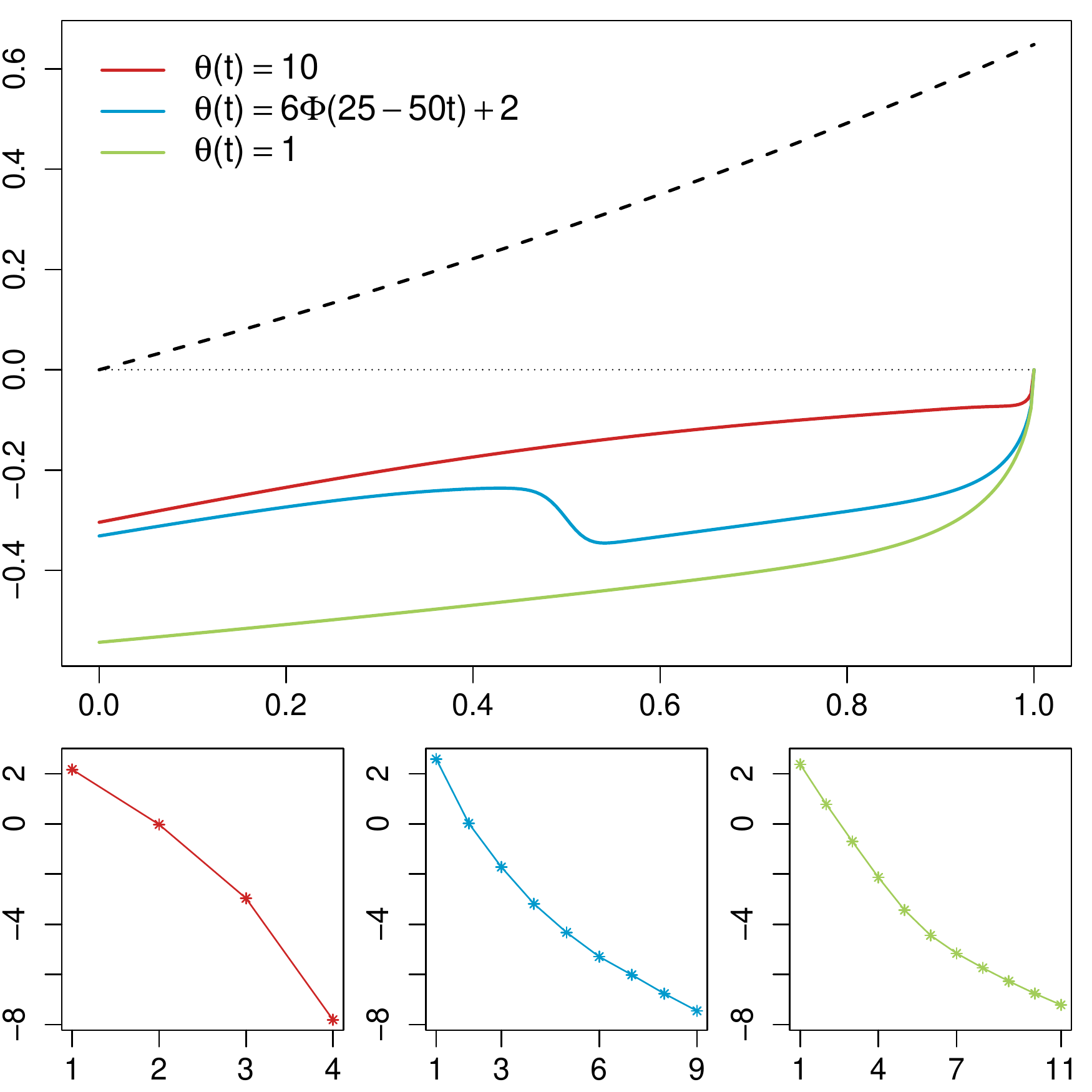}}}\hspace{2pt}
\subfloat[$\theta(t) = 1,\, \alpha(t) = \sin(2\pi t)$.]{%
\resizebox*{0.32\textwidth}{!}{\includegraphics{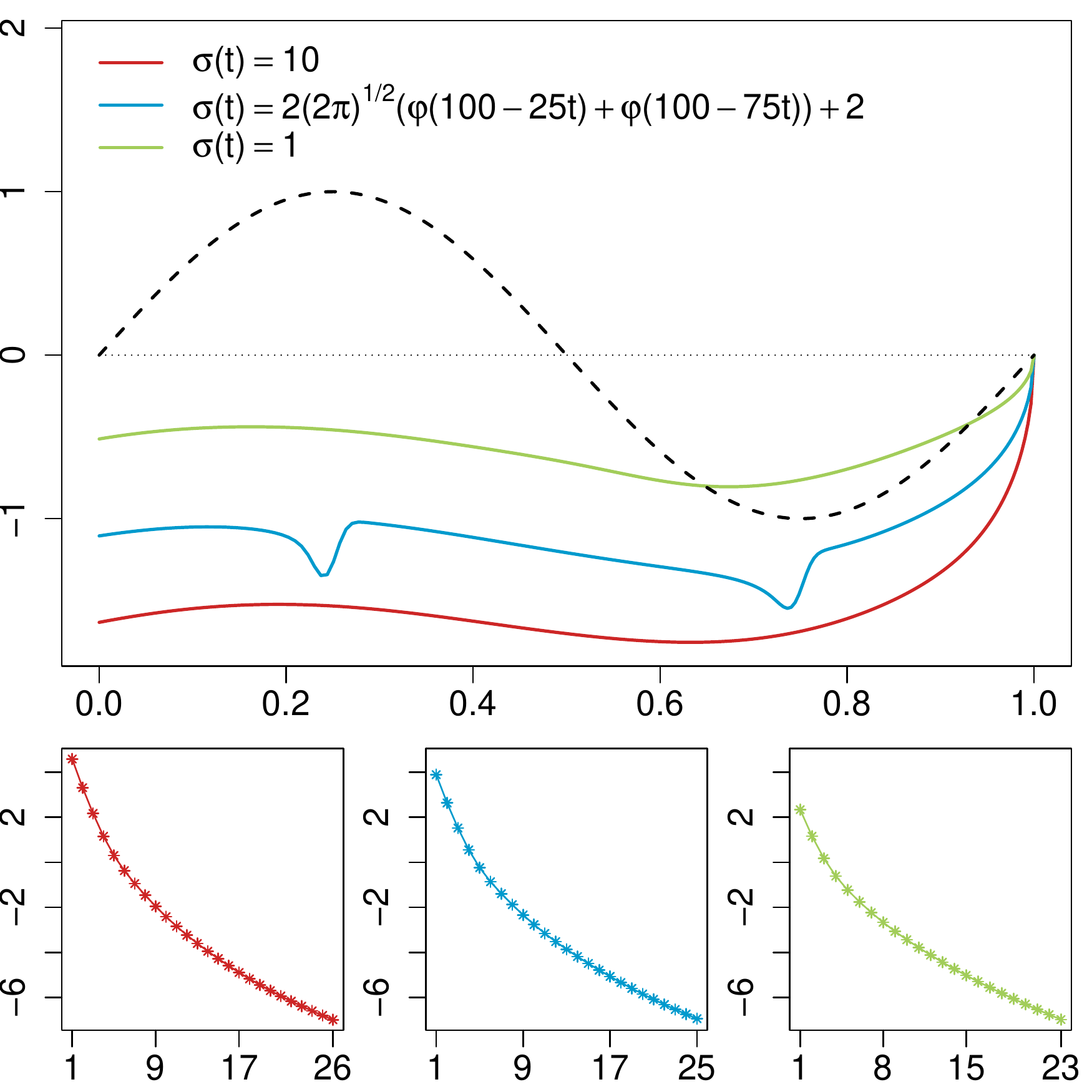}}}
\caption{\small 
For the images on top, the solid colored lines represent the OSBs for the different choices of the $\alpha$, $\theta$, and $\sigma$ functions shown in the legend. The dashed lines represent the pulling level $\alpha$ and the dotted line is placed at the strike price level $A$. $\Phi$ and $\phi$ represent the distribution and density functions of a standard normal. We set $T = 1$, $\lambda = 1$, and $A = 0$. The smaller images below provide the errors $\log_{10}(d_k)$ between consecutive boundaries for each iteration $k = 1, 2,\ldots$ of \eqref{eq:Picard}.} 
\label{fig:OSB_shape}
\end{figure}

Figure \ref{fig:OSB_shape} provides insights into the shape of the OSB for different sets of slopes and pulling levels. We also allow the process to have non-constant volatility, as discussed in Remark \ref{rmk:time_dependent_volatility}. In each image, two of the three functions that define the dynamics in \eqref{eq:GM_SDE} are fixed, and the remaining one varies to make clear the marginal effect it has on the OSB's shape. We use the functions $\Phi$ and $\phi$ to model abrupt changes in a smooth fashion. Figure \ref{fig:OSB_shape}(a) models a regime change of the stock price, where the switching happens smoothly, deterministically, and depends only on the time variable. It also reflects the attracting behavior of the OSB towards $\alpha$. Figure \ref{fig:OSB_shape}(b) exemplifies how this attraction strengthens as $\theta$ increases, and also depicts a sudden increase of the pulling strength, which could represent a sharp increase of confidence in a belief of the price evolution. Figure \ref{fig:OSB_shape}(c) introduces periods where the volatility spikes to unusual levels before returning to a baseline. It also shows how the OSB is repelled from the pulling level as the volatility increases, which, in agreement with Remark \ref{rmk:time_dependent_volatility}, coincides with the effect of reducing the slope.

Recent interest in OSPs related to diffusion bridges whose drift is linear in the space component (see, e.g., \cite{DAuria-2020-discounted} and \cite{DAuria-2021-optimal}) led us to investigate the applicability of our OSB within this context. The difficulty of dealing with these processes is that their slopes $\theta$ explode as the time approaches the horizon, hence they do not fit our assumption of a bounded slope. However, we suggest in Figure \ref{fig:OSB_approximation} that this inconvenience could be circumvented in practice by arbitrarily approximating the explosion with non-exploding smooth drifts of time-dependent OU processes. We do so by relying on the OSB of a Brownian Bridge (BB) and an Ornstein--Uhlenbeck Bridge (OUB) in the non-discounted scenario, whose solutions are available for some particular settings. 

\begin{figure}
\centering
\subfloat[$\alpha(t) = 0,\, \sigma(t) = 1$.]{%
\resizebox*{0.48\textwidth}{!}{\includegraphics{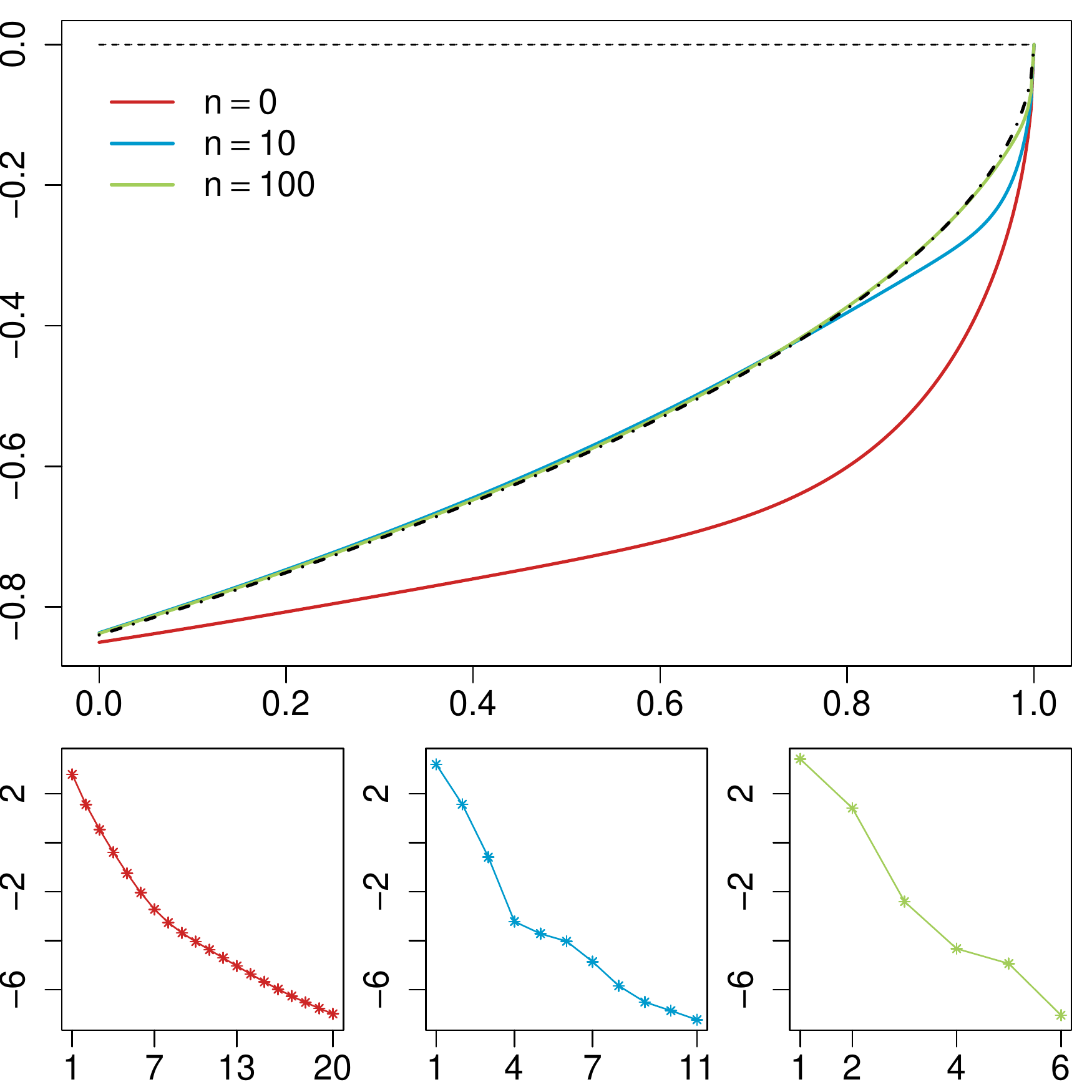}}}\hspace{5pt}
\subfloat[$\alpha(t) = 2\sech(5(1 - t)),\, \sigma(t) = 1$.]{%
\resizebox*{0.48\textwidth}{!}{\includegraphics{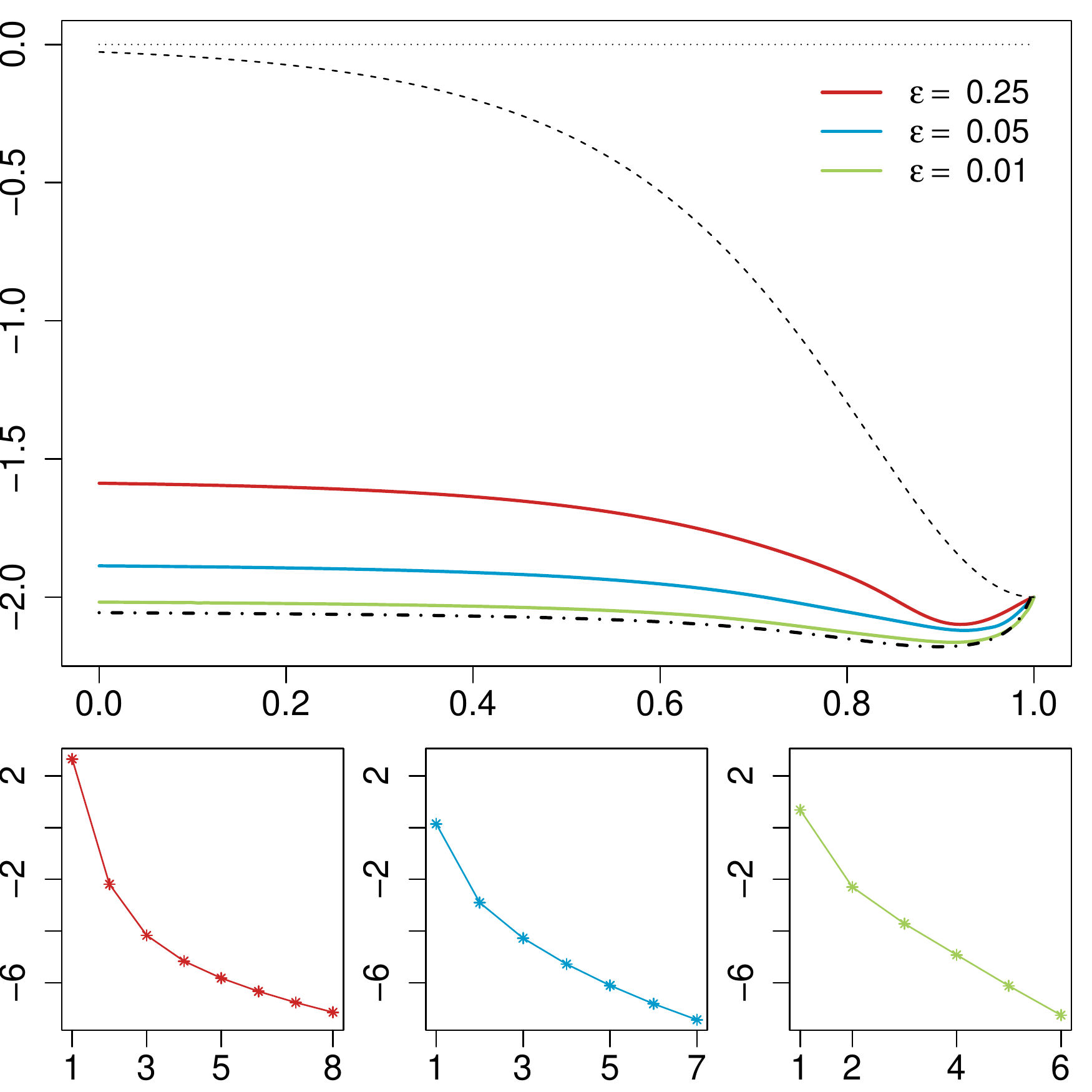}}}
\caption{\small For the images on top, the solid colored lines represent the OSBs for different $n$ and $\varepsilon$ (see main text). The dotted line is placed at the strike price level $A$, while the dashed line stands for the pulling level $\alpha$, and the dotted-dashed lines are the OSBs of the BB (a) and the OUB (b). We set $A = 0$, $\lambda = 0$, and $T = 1$. The smaller images below are analogous to those in Figure \ref{fig:OSB_shape}.}
\label{fig:OSB_approximation}
\end{figure}

The OSB for the BB with pinning point $A$ and gain function $x\mapsto (A - x)^+$ takes the form $A - B\sigma\sqrt{T - t}$, for $B \approx 0.8399$, a result that goes back to the work of \cite{Shepp-1969-explicit} and its connection to the American option by \cite{DAuria-2020-discounted}. We use, in Figure \ref{fig:OSB_approximation}(a), the Taylor-expansion approximation (at $t=0$) of the BB slope $(1 - t)^{-1}$, namely, $\sum_{i=0}^n t^i$, for different values of the order $n$ indicated in the legend. The results in Figure \ref{fig:OSB_approximation}(a) show that the application of the Picard iteration \eqref{eq:Picard} on the time-dependent OU that order-$n$ approximates a BB converges to the expected, known, OSB.

The OSB of the OUB is not known in closed-form, but rather as a characterization via the free-boundary equation provided in \cite{DAuria-2021-optimal}, which can be also solved by means of a Picard algorithm. In that paper the authors treated the OUB case with the identity as the gain function. However, a simple translation with respect to $A$ provides the OSB for $x\mapsto (A - x)^+$ as the gain function whenever the pinning point is lower than $A$. This is easy to note as both gain functions, $x\mapsto A - x$ and $x\mapsto (A - x)^+$ coincide for all $x \leq A$. Then, there is no fundamental difference in both OSPs (for $\lambda = 0$) whenever the boundary associated with the OUB lies below $A$. In Figure \ref{fig:OSB_approximation}(b), to avoid the explosion of the OUB slope, $t\mapsto a\coth(a(1 - t))$, we use instead the smooth function $t\mapsto \Ind(t \leq 1 - \varepsilon)\blrc{a\coth(a(1 - t))} + \Ind(t > 1 - \varepsilon)\blrc{1 - \exp\{-a^2(t - 1 + \varepsilon)/\sinh^2(a\varepsilon)\} + a\coth(a\varepsilon)}$, for $a = 5$, and different values of $\varepsilon$ specified in the legend. Figure \ref{fig:OSB_approximation}(b) shows that the application of \eqref{eq:Picard} on the time-dependent OU that approximates an OUB with pinning point $z = 2$ converges to the (numerically-computed) OSB of the OUB.

Finally, Figure \ref{fig:OSB_convergence_and_discount} illustrates how the algorithm's output seems to converge as the partition size increases. Observe the relatively good performance with just a few points in the partition. Additionally, we show how changing the discounting rate affects the OSB's shape: larger discounts decrease the separation of the OSB with respect to the strike price.

\begin{figure}
\centering
\subfloat[$N = 5$.]{%
\resizebox*{0.32\textwidth}{!}{\includegraphics{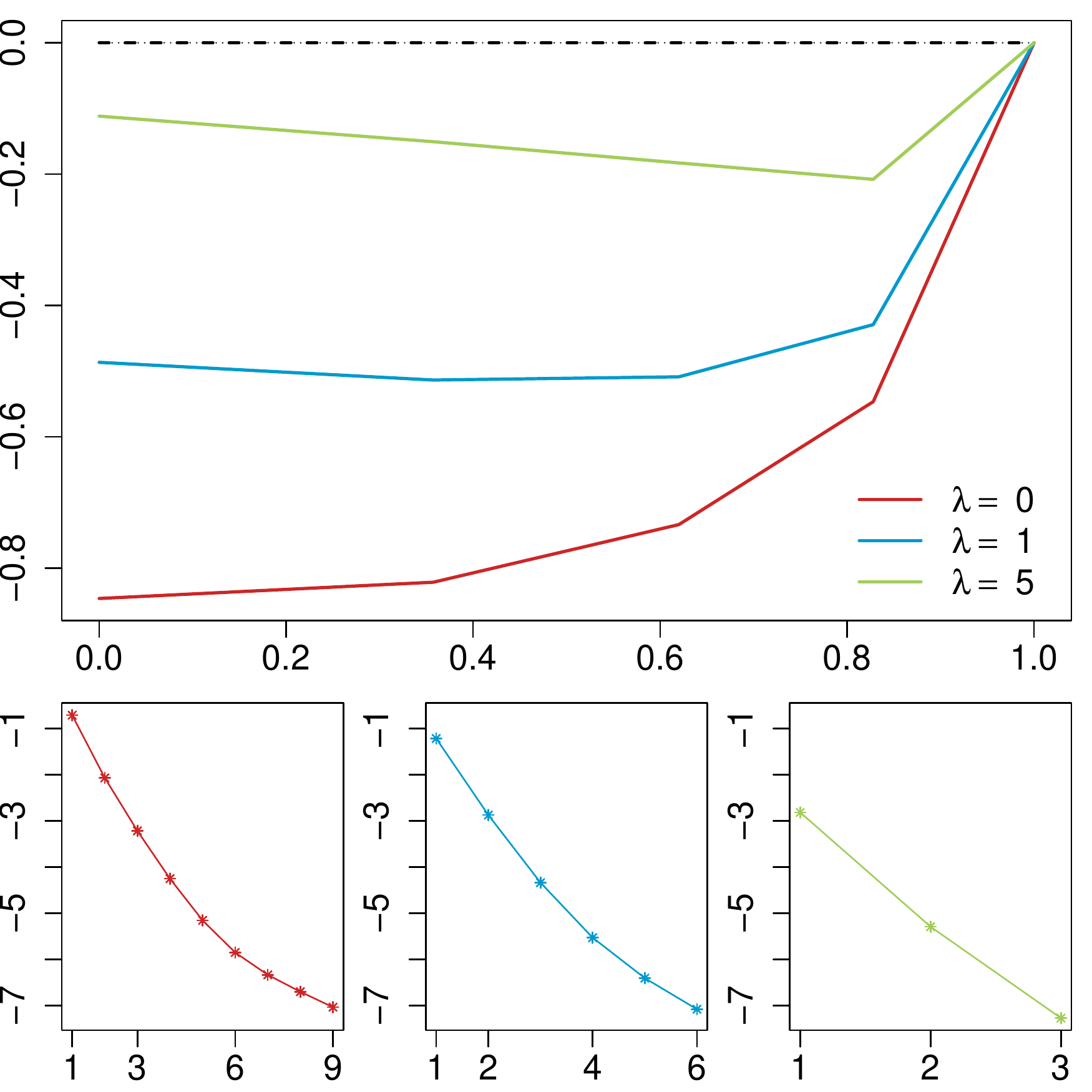}}}\hspace{5pt}
\subfloat[$N = 20$.]{%
\resizebox*{0.32\textwidth}{!}{\includegraphics{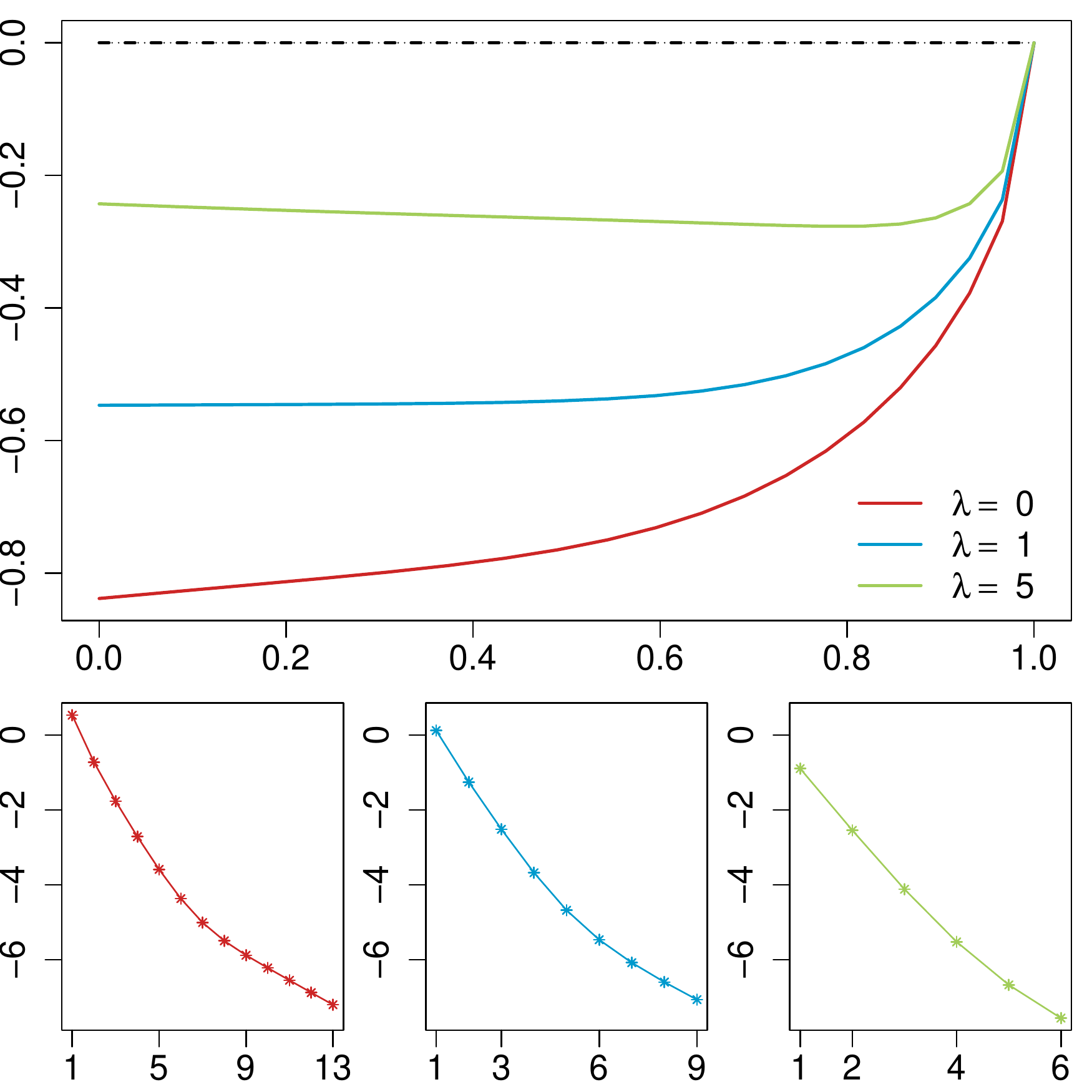}}}
\subfloat[$N = 200$.]{%
\resizebox*{0.32\textwidth}{!}{\includegraphics{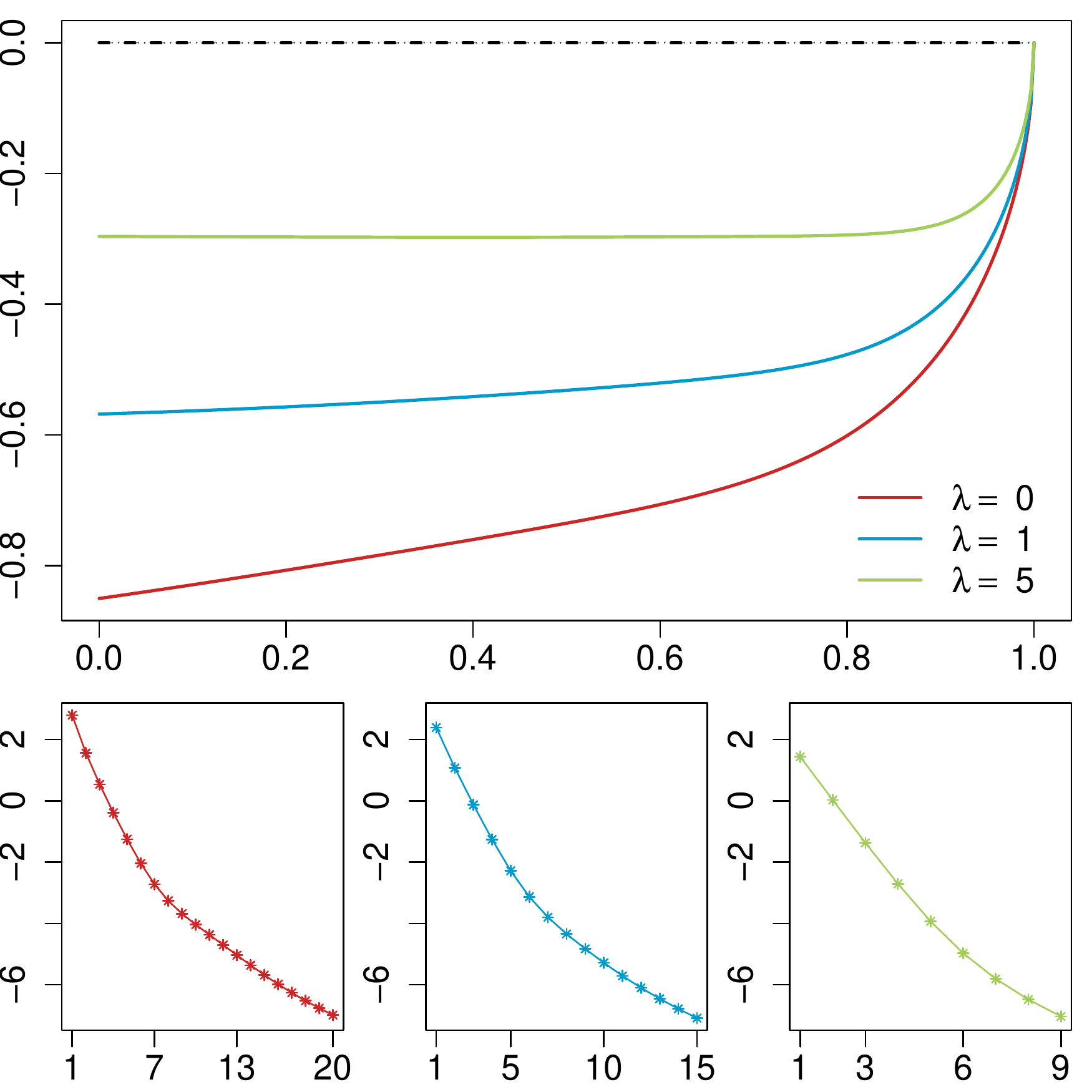}}}
\caption{\small For images on top, the solid colored lines represent the OSBs for different values of the discounting rate specified in the legend, while the dashed line is placed at $A = 0$. We set $\theta(t) = 1$, $\alpha(t) = 0$, $\sigma(t) = 1$, $A = 0$, and $T = 1$. The smaller images below are analogous to those in Figure~\ref{fig:OSB_shape}.} 
	\label{fig:OSB_convergence_and_discount}
\end{figure}

The repository \url{https://github.com/aguazz/AmOpTDOU} contains all the required R scripts to reproduce the numerical experiments.

\section{Concluding remarks}\label{sec:conclusions}

If a time-dependent Ornstein--Uhlenbeck model fits well the price of the underlying asset, the strategy that maximizes the mean gain of a holder of an American put option on that asset is that of exercising once the price lies below the solution of the integral equation in Theorem \ref{th:fb}. This integral equation can be solved numerically using the algorithm from Section \ref{sec:numerical_results}.

The problem of the optimal exercise of the American option was treated in this paper as an optimal stopping problem. We rely on a probabilistic methodology similar to that used in \cite{Peskir-2005-American}. That is, we obtained enough regularity conditions on the value function and on the optimal stopping boundary to apply an extension of the Itô formula and come up with the free-boundary (Volterra integral) equation. Contrary to \cite{Peskir-2005-American} and many papers that followed the same approach (see, e.g., \cite{Peskir-2003-Asian, Peskir-2005-Rusian, Glover-2010-British, Glover-2011-British, Kitapbayev-2014-lookback, DeAngelis-2020-optimal}), our optimal stopping boundary is not necessarily monotonic (see Figure \ref{fig:OSB_shape}), which makes more complicated the derivation of certain properties,  especially the smooth-fit condition.  
To overcome this handicap, we obtained the Lipschitz continuity of the boundary by adapting the work of \cite{DeAngelis-2019-Lipschitz} to fit our non-differentiable gain function and relax other restrictive assumptions. We then proved that Lipschitz continuity suffices to obtain that the underlying process enters the stopping set immediately after starting on the optimal stopping boundary and, then, by relying on the work of \cite{DeAngelis-2020-global}, we derived the smooth-fit condition. It is worth highlighting the comparison method used to get the lower bound of the optimal boundary in Proposition \ref{pr:boundary_existence}, as it seems to be extensible to other settings.

Numerical experiments performed in Section \ref{sec:numerical_results} revealed a wide flexibility in the shape of the optimal boundary generated by changing the form of the time-dependent parameters. It also showed that the errors produced by the Picard iteration algorithm decrease, suggesting that the fixed-point operator in the free-boundary equation \eqref{eq:free-boundary_equation_gaussian} could be a contracting map. 

\section*{Funding}

The authors acknowledge support from grants PID2020-116694GB-I00 (first and second authors), and PGC2018-097284-B-100 (third author), funded by MCIN/AEI/10.13039/\-501100011033 and by ``ERDF A way of making Europe''. The research of the second and third was also supported by the Community of Madrid through the framework of the multi-year agreement with Universidad Carlos III de Madrid in its line of action ``Excelencia para el Profesorado Universitario'' (EPUC3M13). 

\section*{Acknowledgement}

The authors acknowledge the comments of two anonymous referees that led to significant improvements in the manuscript.

\bibliographystyle{apalike-custom}
\bibliography{AmOpOU-arXiv}

\end{document}